 \documentclass{amsart}

\usepackage{epsfig}  		
\usepackage{epic,eepic}       

\usepackage{a4wide,amssymb,latexsym,amsthm,amsmath}
\usepackage{eucal}
\usepackage{graphicx}
\usepackage{algorithm} 
\usepackage{algorithmic}
\usepackage{graphicx}
\usepackage{subcaption}
\usepackage{bbm}
\usepackage{lipsum}

\newcommand{\argmin}{\mathrm{arg\,min}}

\theoremstyle{definition}

\theoremstyle{remark}

\numberwithin{equation}{section}

\begin{document}


\title[A shear-compression damage model for block caving]{A shear-compression damage model for the simulation of underground mining by block caving}


\author[Bonnetier]{Eric Bonnetier}
\address{Fourier Institute, Universit\'e Grenoble-Alpes, 700 Avenue Centrale, 38401 Saint-Martin-d'H\`eres, France.}
\email{eric.bonnetier@univ-grenoble-alpes.fr}

\author[Gaete]{Sergio Gaete}
\address{Codelco Chile, Divisi\'on El Teniente, Mill\'an 1020, Rancagua, Chile.}
\email{sgaete@codelco.cl}

\author[Jofre]{Alejandro Jofre}
\address{Centro de Modelamiento Matem\'atico (CMM) UMI 2807 CNRS-UChile and Departamento de Ingenier\'ia Matem\'atica (DIM), Universidad de Chile, Beauchef 851, Casilla 170-3, Correo 3, Santiago, Chile.}
\email{ajofre@dim.uchile.cl}

\author[Lecaros]{Rodrigo Lecaros}
\address{Departamento de Matem\'atica, Universidad T\'ecnica Federico Santa Mar\'ia, Casilla 110-V, Valpara\'iso, Chile}
\email{rodrigo.lecaros@usm.cl}

\author[Montecinos]{Gino Montecinos}
\address{Departamento de Ciencias Naturales y Tecnolog\'ia,  Universidad de Ays\'en, Obispo Vielmo 62, 5950000 Coyhaique, Chile.}
\email{gino.montecinos@uaysen.cl}

\author[Ortega]{Jaime H. Ortega}
\address{Centro de Modelamiento Matem\'atico (CMM) UMI 2807 CNRS-UChile and Departamento de Ingenier\'ia Matem\'atica (DIM), Universidad de Chile, Beauchef 851, Casilla 170-3, Correo 3, Santiago, Chile.}
\email{jortega@dim.uchile.cl}

\author[Ram\'irez-Ganga]{Javier Ram\'irez-Ganga}
\address{Centro de Modelamiento Matem\'atico (CMM) UMI 2807 CNRS-UChile and Departamento de Ingenier\'ia Matem\'atica (DIM), Universidad de Chile, Beauchef 851, Casilla 170-3, Correo 3, Santiago, Chile.}
\email{jramirez@dim.uchile.cl}

\author[San Mart\'in]{Jorge San Mart\'in}
\address{Centro de Modelamiento Matem\'atico (CMM) UMI 2807 CNRS-UChile and Departamento de Ingenier\'ia Matem\'atica (DIM), Universidad de Chile, Beauchef 851, Casilla 170-3, Correo 3, Santiago, Chile.}
\email{jorge@dim.uchile.cl}



\begin{abstract}
Block caving is an ore extraction technique used in the copper mines of Chile.
It uses gravity to ease the breaking of rocks, and to facilitate the extraction
from the mine of the resulting mixture of ore and waste.
To simulate this extraction process numerically and better understand its impact 
on the mine environment, we study 3 variational models for damage, based on the gradient damage model of Pham and Marigo \cite{pham2010approche1,pham2010approche2}. In these models, the damage criterion may exhibit an anisotropic dependence on the spherical and deviatoric parts of the stress tensor. We report simulations that satisfactorily represent the expected evolution of the stress field in a block caving operation.
\end{abstract}


 \maketitle



\section{Introduction} 
Block caving is a mining method in which ore blocks are undermined, causing the rocks 
to cave, thus allowing broken ore to be removed with greater ease at draw-points (see \cite{hartman2002introductory}).  This method is based on the sinking principle and mineral breakage due to the removal of a large underground supporting area of rock.  The extraction of material from this cavity induces a reconfiguration  of the stress field and the fracking of the mineral above it, which then, under the action of gravity, falls into the generated vacuum generated.
The friction and resulting attrition during the fall, defines the
size of the blocks of mineral at the extraction point. From there, the material is transferred to the processing plant for further size reduction or treatment.

In this paper we report the numerical study of a mathematical model, 
which describes the effect of such a block caving process in the rock mass.
The fracking of the rocks is considered from the point of view of a continuum
approach, in which the cracks that weaken the material stiffness form 
a damaged zone. Specifically, we seek to recover the evolution of such damage in the rock mass around the cavity, so as to understand the subsidence of the cavity ceiling \cite{brown2002block}.

The theory of fracture rests on the assumption that the growth of a crack
requires an amount of work proportional to the newly created surface area.
The corresponding proportionality constant is known as fracture toughness.
As the crack grows, the displacement field is assumed to instantly 
find in a new equilibrium configuration. Since the displacement may be 
discontinuous across the crack increment, the stored elastic energy 
decreases as the crack grows. Classical models of fracture are based on the Griffith criterion, which stipulates that the energy release rate  is equal to the work required to create the crack increment, i.e.
to the fracture toughness. From a mathematical point of view, the associated theory of evolution is not satisfactory, as it assumes before-hand knowledge of the crack path, and can
only be made rigorous for smooth crack topologies.

The variational approach of G. Francfort and J.J. Marigo \cite{francfort1998revisiting,bourdin2008variational} has brought a paradigm shift to the theory of (quasistatic) brittle fracture, which has allowed to lift some of these mathematical hurdles. In their work, the elastic displacement is defined via the minimization of a potential energy, which is the sum of the stored elastic energy and the surface energy  of sets of discontinuity, the latter defining the fracture path. The potential energy functional takes the general form
\begin{equation}
	\mathcal{P}(u) =
	\frac{\mu}{2} \int_\Omega \vert\nabla u\vert^2 dx + G_c \mathcal{H}^{n-1}(S(u)),
\end{equation}
for displacements $u \in SBV(\Omega)$ the special functions of bounded variations.
Here $\Omega$ represents the set initially occupied by the elastic body,  $\mu$ is a constant associated to the material stiffness, and $G_c$ denotes the fracture toughness and $S(u)$
the discontinuity set of u. The surface energy is expressed in terms of $n-1$ dimensional Hausdorff measure of the discontinuity set $\mathcal{H}^{n-1}(S(u))$.  The evolution of quasi-static fracture in this form is governed by  a condition of stationarity, a condition of irreversibility and an 
energy balance.  Existence and well-posedness results, in particular in relation with the choice of admissible displacement fields, can be found in \cite{francfort1998revisiting,bourdin2008variational,ambrosio1997energies,dal2002model,francfort2003existence}.

In view of the numerical approximation of the above energy, and 
inspired by the Mumford-Shah functional for image segmentation, 
one may introduce an additional variable $v$ whose role is to smear 
the discontinuities of the displacement and track their position.
The Gamma-convergence results of Ambrosio and Tortorelli \cite{ambrosio1990approximation} and the results of \cite{giacomini2005ambrosio} show that an energy of the form
\begin{equation} \label{model_isotropic}
\mathcal{P}_{\varepsilon}(u,v) = \frac{\mu}{2}\int_\Omega (v^2 + \eta_{\varepsilon}) \vert \nabla u\vert^2 dx
+ G_c \int_\Omega \left( \vert\nabla u\vert^2 + \frac{(1-v)^2}{4\varepsilon} \right) dx,
\end{equation}
$\Gamma$-converges to the Griffith energy $\mathcal{P}(u)$.
This approximation is at the basis of the numerical implementation
for the quasistatic evolution of fracture, see \cite{bourdin2000numerical,bourdin2007numerical,bourdin2008variational}.

In a gradient damage model, the failure is represented by an internal variable which alters the stiffness of the material. The term `gradient' refers to the addition of a regularizing term in the internal energy, that depends on the gradient of this internal variable and on an internal length, which fixes the scale of the regularization \cite{pham2010approche1,pham2010approche2,lancioni2009variational,marigo2016overview}. As such, a gradient damage model of evolution may be understood as an
elliptic approximation to a quasistatic fracture evolution:
As the internal length tends to zero, the minima of the gradient damage
energy functional converge (in the proper sense) to minima of the 
Griffith energy functional of brittle fracture \cite{braides1998approximation}.

The aim of this work is to study whether a gradient damage model can 
account for the subsidence of material above an extraction cavity in an 
underground mine. Indeed, one would like to understand the impact of creating an extraction cavity on the state of the material above it, so as to assess
security issues for the mining operation or so as to optimize extraction.
Our choice of a regularized variational model is motivated by the fact that
such formulation allows tracking of the damage 
without requiring explicit description of the geometry of the damaged zone.
We hope that the output of such model may be used as indications of the location
of the weak areas in underground mines. It may also prove helpful
as a means to extract information on the state of the material above a cavity 
from onsite measurements.

Our first attempt with the simplest isotropic model of the form \cite{pham2011gradient,marigo2016overview} leads to results that are not very satisfactory. As compression is larger below than above the cavity, most of the damage produced will be located below it, which would not account for a block caving experiment. We therefore introduce a model in which damage affects resistance to shear and compression in a different manner, and we report numerical results that show that damage can occur where it is expected, s0 that such model would be pertinent for the modeling of block caving.

This paper is organized as follows. In Section \ref{DModel}, variational gradient damage models are introduced. We consider 3 choices of damage criterion (isotropic, shear and shear-compression). We describe the numerical method that we use to simulate the evolution of these models in section \ref{NumSec} and report the 2D test cases that we performed. We show how the shear-compression model is more relevant to the modeling of block-caving. Finally, Section \ref{Con} contains a few words of conclusion and comments.

In all the work we will use the summation convention on repeated indices: Vectors and second order tensors are indicated by a lowercase letter, such as $u$ and $\sigma$ for the displacement field and the stress field. Their components are denoted by $u_i$ and $\sigma_{ij}$. Third or fourth order tensors as well as their components are indicated by a capital letter, such as $A$ or $A_{ijkl}$ for the stiffness tensor. Such tensors are considered to be linear maps acting on vectors or second order tensors and this action is denoted without dots, such as $A\varepsilon$ whose $ij$-component is $A_{ijkl}\varepsilon_{kl}$. The inner product between two vectors or two tensors of the same order is indicated by a colon. For instance $a:b$ stands for $a_ib_i$ and $\sigma : \varepsilon$ for $\sigma_{ij}\varepsilon_{ij}$. We use the notation $A>0$ to denote a positive definite tensor. The time derivative is indicated by a dot, for instance $\dot{\alpha}:=\frac{\partial \alpha}{\partial t}$.

\section{Damage models for block caving process}\label{DModel}

\subsection{General presentation of damage models}

We briefly explain the construction of the gradient damage models, which inspired our models for block caving. We refer the readers to \cite{marigo2016overview} and references therein for a thorough review of its variational and constitutive ingredients as well as its properties especially when applied to brittle fracture in a quasi-static setting. 

Consider a homogeneous $n$-dimensional body whose reference configuration is 
the open connected bounded set $\Omega \subseteq \mathbb{R}^n$. The vector field $u$ denotes the displacement and damage is represented by a scalar parameter $\alpha \in [0,1]$. The state of a volume element is characterized by the triplet 
$(\varepsilon,\alpha,\nabla \alpha)$, where  $\varepsilon(u)=\frac{1}{2}\left(\nabla u + \nabla u^T \right)$ denotes the linearized strain tensor. When $\alpha=0$, the material is in the healthy-undamaged state, and behaves like an isotropic linear elastic material. When $\alpha=1$, the material is in a fully damage state. Its behavior is characterized by a state function 
\begin{eqnarray}\label{statefunc}
W_{\ell}:(\varepsilon,\alpha,\nabla \alpha) &\mapsto & 
\psi(\varepsilon,\alpha) + w(\alpha) + \frac{1}{2}w_1\ell^2\vert \nabla \alpha \vert^2.
\end{eqnarray}
which is the sum of three terms: the stored elastic energy $\psi(\varepsilon,\alpha)$, the local part of the dissipated energy by damage $w(\alpha)$ and the non local part of the dissipated energy.
The stored elastic energy takes the form
\begin{equation}
\psi(\varepsilon,\alpha)=\frac{1}{2}A(\alpha)\varepsilon(u):\varepsilon(u),
\end{equation}
where $A(\alpha)$ represents the rigidity tensor of the material, given the state of damage $\alpha$. The material becomes less rigid as $\alpha$ increases, that is, 
\begin{equation}
A(\alpha)=\mathrm{a}(\alpha) A_0,
\end{equation}
where $\alpha \mapsto \mathrm{a}(\alpha)$ is a decreasing scalar function
such that $\mathrm{a}(0)=1$ and $\mathrm{a}(1) = 0$, and where $A_0$ is the isotropic elasticity tensor with Young's modulus $E$ and Poisson ratio $\nu$. The local dissipated energy density is a positive, increasing function of $\alpha$, with $w(0) = 0$ and $w(1) = w_1$. The value $w_1$ thus represents the energy dissipated during a complete, homogeneous damage process ($\nabla \alpha =0$) of a volume element.

For simplicity, the non local dissipated energy density is assumed to be a 
quadratic function of the gradient of damage. The parameter $\ell > 0$ can be considered as an internal length characteristic of the material, which controls the thickness of the damage localization zones.

\subsection{The evolution problem}

We assume that the body $\Omega$ is subjected to a time dependent loading $U(\cdot,t)$, $F(\cdot,t)$ and $f(\cdot,t)$ which consists in an imposed displacement on the part of the boundary $\partial \Omega_U$, surface forces on the complementary part of the boundary $\partial \Omega_F$, and volume forces, $t$ denoting the time parameter. The space of kinematically admissible displacement fields at time $t$ is the set
\begin{equation}
\mathcal{C}_t=\left\{ v \in (H^1(\Omega))^n \mbox{ : } v=U(\cdot,t) \mbox{ on } \partial \Omega_U \right\},
\end{equation}
while we seek the damage field in the set $\mathcal{D}$, where
\begin{equation} \label{defD}
\mathcal{D}= \left\{ \beta\in H^1(\Omega) \mbox{ : } \beta(x) \in [0,1] \mbox{ for almost all } x \in \Omega \right\}.
\end{equation}

The law of evolution of the damage in the body is written in variational form.
If $(u,\alpha)$ denotes a pair of admissible displacement and damage fields at time $t$, i.e., if $u \in \mathcal{C}_t$ and $\alpha \in \mathcal{D}$, then the total energy of the body at time $t$ in this state is given by
\begin{equation}\label{totalenergy}
\mathcal{P}_t(u,\alpha)= \int_{\Omega} W_{\ell}(\varepsilon(u(x,t)),\alpha(x,t),\nabla \alpha(x,t))dx- \int_{\Omega} f(x,t) \cdot u(x,t) dx - \int_{\partial \Omega_F} F(x,t) \cdot u(x,t) dS. 
\end{equation}

In the variational approach, the quasi-static evolution for the displacement and the damage is formulated as a first-order unilateral minimality condition on the functional \eqref{totalenergy}, under the condition of irreversibility that damage can only increase \cite{marigo1989constitutive}, i.e., that $\dot{\alpha}\geq 0$. The evolution is stated in the following form: Find $(u,\alpha)\in \mathcal{C}_t\times \mathcal{D}$ such that
\begin{equation}\label{firstordercond}
\forall\; (v,\beta) \in \mathcal{C}_t \times\mathcal{D}_t(\alpha), \quad\quad
\mathcal{P}'_t(u,\alpha)(v-u,\beta-\alpha)\geq 0,
\end{equation}
where $\mathcal{D}_t(\alpha)=\left\{ \beta \mbox{ : } \alpha(x,t) \leq \beta(x) \leq 1 \mbox{ in } \Omega\right\}$ and $\mathcal{P}'_t(u,\alpha)(v,\beta)$ denotes the directional derivative of $\mathcal{P}_t$ at $(u(\cdot,t),\alpha(\cdot,t))$ in the direction $(v,\beta)$, i.e., the linear form defined by
\begin{equation}
\begin{split}
\mathcal{P}'_t(u,\alpha)(v,\beta) =& \int_{\Omega}  \mathrm{a}(\alpha)A_0\varepsilon(u)  : \varepsilon(v) dx \\ &  +  \int_{\Omega} \left(  \left(\frac{1}{2}\mathrm{a}'(\alpha)A_0\varepsilon(u):\varepsilon(u) + w'(\alpha)\right) \beta  +w_1\ell^2 \nabla \alpha \cdot \nabla \beta \right) dx \\& - \int_{\Omega}f\cdot v-\int_{\partial \Omega_F}F \cdot v dS.
\end{split}
\end{equation}

Assuming existence and regularity of a solution, a local formulation can be deduced from this global variational formulation by integration by parts and classical localization arguments.

The equilibrium equation  is obtained by testing the variational inequality (\ref{firstordercond}) for $\beta=\alpha(\cdot,t)$ and $v\in \mathcal{C}_t$. This yields the system of quasi-static equilibrium equations
\begin{equation}
\begin{array}{r c l c l}
\mathrm{div}(A(\alpha)\varepsilon(u)) +f  &= &0 & \mbox{ in } &\Omega, 
\\ [8pt]
A(\alpha)\varepsilon(u) \cdot n &= &F &\mbox{ on } &\partial \Omega_F, 
\\ [8pt]
u &= &U &\mbox{ on } &\partial \Omega_U.
\end{array}
\end{equation}

In addition, an equation for $\alpha$ is obtained by testing (\ref{firstordercond}) for arbitrary $\beta$ in the convex cone $\mathcal{D}_t(\alpha)$ with $v = u(\cdot,t)$, leading to the variational inequality governing the evolution of the damage
\begin{eqnarray}\label{damagecriterion}
\forall \beta \in \mathcal{D}_t(\alpha), & \displaystyle \int_{\Omega} \left( \left(\frac{1}{2}\mathrm{a'}(\alpha)A_0\varepsilon(u):\varepsilon(u) + w'(\alpha)\right)(\beta-\alpha) +  w_1\ell^2 \nabla \alpha \cdot \nabla (\beta-\alpha)\right)dx \geq 0.
\end{eqnarray}

After integration by parts, we find the strong formulation for the damage evolution problem in the form of the Kuhn-Tucker conditions for unilaterally constrained variational problems
\begin{enumerate}
	\item {\bf Irreversibility:} $\dot{\alpha} \geq 0$ in $\Omega$.
	\item {\bf Damage criterion:} 
	$\displaystyle\frac{1}{2}\mathrm{a}'(\alpha)A_0\varepsilon(u):\varepsilon(u)
	+ w'(\alpha)- w_1\ell^2 \Delta \alpha \geq 0$ in $\Omega$.
	\item {\bf Energy balance:} 
	$ \displaystyle \dot{\alpha}\left( \frac{1}{2}\mathrm{a}'(\alpha)
	A_0\varepsilon(u):\varepsilon(u)
	+ w'(\alpha)- w_1\ell^2 \Delta \alpha \right) =0$ in $\Omega$.
	\item {\bf Boundary Conditions:} $ \displaystyle  \frac{\partial\alpha}{\partial n} 
	\geq 0$ and $  \displaystyle \dot{\alpha} \frac{\partial\alpha}{\partial n}=0$ on $\partial\Omega$.
\end{enumerate}

The energy balance states that, at each point, damage only increases if the damage yield criterion is satisfied, that is if the damage criterion is an equality.

\subsection{Generalized standard materials}\label{GSM}

In this section, we cast the damage models we studied in the framework of generalized standard materials \cite{halphen1975materiaux}. In an isothermal process, the Clausius-Duhem inequality that expresses the second law of thermodynamics takes the form
\begin{equation} \label{ineq_cd}
\Phi = \sigma:\dot{\varepsilon}  - \dot{W} \geq 0,
\end{equation}
where $\Phi$ denotes the density of dissipated power, $W$ the free energy,
$\sigma$ the stress tensor and $\dot{\varepsilon}$ the strain rate tensor,
defined by
\begin{equation}
\dot{\varepsilon} = \varepsilon(\dot{u}) = \frac{1}{2}\left(\nabla \dot{u} + \nabla \dot{u}^T \right).
\end{equation}
In our case, the free energy is assumed to depend on the strain and on an internal variable $\alpha$ that measures damage, \eqref{ineq_cd} can be rewritten in the form
\begin{equation}
\left(\sigma - \displaystyle\frac{\partial W}{\partial \varepsilon}\right):\dot{\varepsilon} - \displaystyle\frac{\partial W}{\partial \alpha}:\dot{\alpha} \geq 0.
\end{equation}

In the theory of generalized standard materials, it is postulated that there exists a non-negative, convex, lower semi-continuous function $\varphi(\dot{\varepsilon},\dot{\alpha})$, which satisfies $\varphi(0,0) = 0$, and such that $\Phi = \dot{\varphi}$. This assumption yields the relations
\begin{eqnarray}\label{sigmaeqdis}
\frac{\partial \varphi}{\partial \dot{\varepsilon}}&=& \sigma  - \displaystyle\frac{\partial W}{\partial \varepsilon},\\ 
\label{alphaeqdis}
\frac{\partial \varphi}{\partial \dot{\alpha}} &=& \displaystyle-\frac{\partial W}{\partial \alpha}.
\end{eqnarray}

We next consider two choices for the free energy and the associated damage models.

\subsubsection{Gradient damage model}\label{GDM}

In a first model, we assume that the free energy and the pseudopotential of dissipation take the following form
\begin{eqnarray}
W(\varepsilon,\alpha) &=& 
\frac{1}{2} \mathrm{a}(\alpha)A_0 \varepsilon(u): \varepsilon(u)+ w(\alpha) 
+ \frac{1}{2} w_1\ell^2 \vert \nabla\alpha \vert^2, \\
\varphi(\dot{\varepsilon}, \dot{\alpha}) &=& I_+(\dot{\alpha}),
\end{eqnarray}
where, $I_+(\dot{\alpha})$ is the indicator function defined by
\begin{equation}
I_+(x) = \left\{
\begin{array}{rcl}
0,      & \mathrm{if } & x \geq 0, \\
+\infty, & \mathrm{if } & x < 0.
\end{array}
\right.
\end{equation}

Equations \eqref{sigmaeqdis}-\eqref{alphaeqdis} then yield
\begin{eqnarray}
\sigma(u,\alpha) - \mathrm{a}(\alpha)A_0\varepsilon & = & 0,\\
\frac{1}{2}\mathrm{a}'(\alpha)A_0 \varepsilon : \varepsilon +  w'(\alpha) - w_1\ell^2 \Delta\alpha & = & - R(\dot{\alpha}),
\end{eqnarray}
where, $R(\dot{\alpha}) = 0$ if $\dot{\alpha}\geq 0$ and $R(\dot{\alpha})\in (-\infty,0]$ if $\dot{\alpha}= 0$.

This gradient damage model is studied in \cite{marigo2016overview} and its evolution is defined as follows:
\begin{enumerate} 
	\item The stress tensor $\sigma(x,t)=\sigma(u(x,t),\alpha(x,t))=\mathrm{a}(\alpha(x,t))A_0\varepsilon(u(x,t))$ satisfies the equilibrium equations
	\begin{equation}
	\begin{array}{r c l c l}
	\mathrm{div}(\sigma(x,t)) +f(x,t)  &= &0 & \mbox{ in } &\Omega, \\
	\sigma(x,t) \cdot n &= &F(x,t) &\mbox{ on } &\partial \Omega_F, \\
	u(x,t) &= &U(x,t) &\mbox{ on } &\partial \Omega_U.
	\end{array}
	\end{equation}
	
	\item The damage field $\alpha(x,t)$ satisfies the nonlocal damage criterion
	\begin{eqnarray}
	\frac{1}{2}\left( \mathrm{a}'(\alpha(x,t))A_0\varepsilon(u(x,t)):\varepsilon(u(x,t))\right)+ w'(\alpha(x,t))-w_1\ell^2 \Delta \alpha(x,t) \geq 0, & \mbox{ in } \Omega,
	\end{eqnarray}
	and the nonlocal consistency condition
	\begin{eqnarray}
	\left(\frac{1}{2}\left( \mathrm{a}'(\alpha)A_0\varepsilon(u(x,t)):\varepsilon(u(x,t))\right)+ w'(\alpha(x,t))-w_1\ell^2 \Delta \alpha(x,t) \right) \dot{\alpha}(x,t) = 0, & \mbox{ in } \Omega.
	\end{eqnarray}	
\end{enumerate}

\subsubsection{Gradient damage model for shear fracture}\label{SDM}

In this model, the pseudopotential of dissipation again takes the form
\begin{equation}
\varphi(\dot{\varepsilon}, \dot{\alpha}) = I_+(\dot{\alpha}).
\end{equation}

For the free energy, we consider the approach based on the orthogonal decomposition of the linearized strain tensor in its spherical and deviatoric components
\begin{eqnarray}
\varepsilon = \varepsilon^s + \varepsilon^d, & \varepsilon^s = \frac{1}{n} \mathrm{tr}(\varepsilon)I, & \varepsilon^d = \varepsilon-\varepsilon^s,
\end{eqnarray}
with $I$ the $n$-dimensional identity tensor. With this decomposition, the free energy $W$ is written as the sum of the spherical and deviatoric contribution of the strain tensor. Moreover, to reproduce the shear fracture, the spherical part remains unaffected by the value of the damage scalar $\alpha$, that is,
\begin{equation}
W(\varepsilon,\alpha) = \left(\lambda + \frac{2 \mu}{n} \right) \frac{\mathrm{tr}(\varepsilon(u))^2}{2}+\mathrm{a}(\alpha)\mu \varepsilon^d:\varepsilon^d + w(\alpha) + \frac{1}{2} w_1\ell^2 \vert \nabla\alpha \vert^2, 
\end{equation} 
where, $\lambda$ and $\mu$ are the Lam\'e coefficients. Then, considering the above mentioned, the equations \eqref{sigmaeqdis}-\eqref{alphaeqdis} yield
\begin{eqnarray}
\sigma(u,\alpha) - \left((2\mu+n\lambda)\varepsilon^s+2\mathrm{a}(\alpha)\mu\varepsilon^d\right) & = & 0,\\
\mathrm{a}'(\alpha)\mu \varepsilon^d : \varepsilon^d +  w'(\alpha) - w_1\ell^2 \Delta\alpha & = & - R(\dot{\alpha}).
\end{eqnarray}    
In this model, proposed by Lancioni and Royer-Carfagni in \cite{lancioni2009variational}, the creation of surface energy may be compensated exclusively by a reduction of the deviatoric elastic energy.

Accordingly, the evolution is defined as follows: at each time $t$,
\begin{enumerate} 
	\item The stress tensor $\sigma(x,t)=\sigma(u(x,t),\alpha(x,t))=\left((2\mu+n\lambda)\varepsilon^s(x,t)+2\mathrm{a}(\alpha(x,t))\mu\varepsilon^d(x,t)\right)$ satisfies the equilibrium equations
	\begin{equation}
	\begin{array}{r c l c l}
	\mathrm{div}(\sigma(x,t)) +f(x,t)  &= &0 & \mbox{ in } &\Omega, \\
	\sigma(x,t) \cdot n &= &F(x,t) &\mbox{ on } &\partial \Omega_F, \\
	u(x,t) &= &U(x,t) &\mbox{ on } &\partial \Omega_U.
	\end{array}
	\end{equation}
	
	\item The damage field $\alpha(x,t)$ satisfies the nonlocal damage criterion
	\begin{eqnarray}
	\mathrm{a}'(\alpha(x,t))\mu \varepsilon^d(x,t) : \varepsilon^d(x,t) + w'(\alpha(x,t))-w_1\ell^2 \Delta \alpha(x,t) \geq 0, & \mbox{ in } \Omega,
	\end{eqnarray}
	and the nonlocal consistency condition
	\begin{eqnarray}
	\left(\mathrm{a}'(\alpha(x,t))\mu \varepsilon^d(x,t) : \varepsilon^d(x,t)+ w'(\alpha(x,t))-w_1\ell^2 \Delta \alpha(x,t) \right) \dot{\alpha(x,t)} = 0, & \mbox{ in } \Omega.
	\end{eqnarray}	
\end{enumerate}

\subsection{A shear-compression damage model}\label{SCDM}

As in Section \ref{GSM}, we assume that the stress-strain relation is given by
\begin{equation}
\sigma(u,\alpha)=\mathrm{a}(\alpha)A_0\varepsilon(u).
\end{equation}
For the equation that governs the evolution of damage, we consider the following relationship
\begin{equation}
\frac{1}{2}H(\varepsilon,\alpha)+w'(\alpha)-w_1\ell^2\Delta \alpha = -R(\dot{\alpha}),
\end{equation}
where $H(\varepsilon,\alpha)$ is a source term for damage.

We define the expression $H(\varepsilon, \alpha)$ in terms of the principal stress of the stress tensor. For a given stress tensor $\sigma=\sigma(u,\alpha)=\mathrm{a}(\alpha)A_0\varepsilon(u)$, we can solve the characteristic equation to explicitly determine the principal values and directions. In three-dimensional case, following \cite{Smith:1961:EST:355578.366316}, if we consider
\begin{eqnarray}
\displaystyle  m  =  \frac{1}{3}\mathrm{tr}(\sigma), &
\displaystyle  q  =  \frac{1}{2}\mathrm{det}(\sigma-mI), &
\displaystyle  p  =  \frac{1}{6} \sum_{ij}(\sigma-mI)^2_{ij},
\end{eqnarray} 
then, from Cardano's trigonometric solutions of $\mathrm{det}[(\sigma - mI)-\lambda I]$ as a cubic polynomial in $\lambda$, the eigenvalues of $\sigma$ are defined by
\begin{eqnarray}
\lambda_1  =  m+2\sqrt{p} \mathrm{cos}(\theta), &
\lambda_2  =  m-\sqrt{p} \left( \mathrm{cos}(\theta) + \sqrt{3} \mathrm{sin}(\theta)  \right), &
\lambda_3  =  m-\sqrt{p} \left( \mathrm{cos}(\theta) - \sqrt{3} \mathrm{sin}(\theta)  \right),
\end{eqnarray}
where $\theta=\frac{1}{3}\mathrm{tan}^{-1}\left( \frac{\sqrt{p^3 - q^2}}{q} \right)$ and $0 \leq \theta \leq \pi$. Thus, if we rewrite the eigenvalues, we have
\begin{eqnarray}
\lambda_1  =  m+2\sqrt{p} \mathrm{cos}(\theta), &
\lambda_2  =  m-2\sqrt{p} \mathrm{cos}\left(\theta-\frac{\pi}{3}\right) , &
\lambda_3  =  m-2\sqrt{p} \mathrm{cos}\left(\theta+\frac{\pi}{3}\right),
\end{eqnarray}
subtracting the eigenvalues, we obtain
\begin{eqnarray}
\lambda_1 - \lambda_2  =  2\sqrt{3p} \mathrm{cos}\left(\theta-\frac{\pi}{6}\right), &
\lambda_1 - \lambda_3  =  2\sqrt{3p} \mathrm{cos}\left(\theta+\frac{\pi}{6}\right) , &
\lambda_2 - \lambda_3  =  2\sqrt{3p}\mathrm{sin}(\theta).
\end{eqnarray}

Then, considering the trigonometric part of the previous equations, we can assert, following the Mohr's circle of stress, that the largest shear component is obtained when $\frac{\lambda_1- \lambda_2}{2} = \sqrt{3p}$.

We assume that the material gets damaged when the magnitude of the shear stress component $S$ is greater than the factor of the magnitude of the normal stress component $N$, i.e., $S > \kappa N$. As mentioned earlier, the maximun shear is $\sqrt{3p}$, our damage criterion will be $\sqrt{3p} > \kappa N$. On the other hand, considering $\sigma = \sigma^s+\sigma^d$, where $\sigma^s$ and $\sigma^d$ are the spherical and deviatoric part of $\sigma$ respectively defined by $\sigma^s=mI$ and $\sigma^d=\sigma-mI$, we have  
\begin{eqnarray}
\vert m \vert = \sqrt{m^2} = \sqrt{\frac{1}{3} \sigma^s : \sigma^s} , &
\sqrt{3p} = \sqrt{\frac{1}{2} \sum_{ij}(\sigma-mI)^2_{ij}} = \sqrt{\frac{1}{2} \sigma^d : \sigma^d},
\end{eqnarray} 
so that damage takes place when
\begin{equation}\label{conddamaged}
\sigma^d:\sigma^d - \frac{2}{3} \kappa \sigma^s:\sigma^s > 0.
\end{equation}

In order to compare our model with the models presented in Section \ref{GSM}, we rewrite the above quantity using the stress-strain relation $\sigma=\mathrm{a}(\alpha)\left(2 \mu \varepsilon +\lambda \mathrm{tr}(\varepsilon)I\right)$, which yields
\begin{eqnarray}
\sigma^d  =  \mathrm{a}(\alpha)2 \mu \varepsilon^d, &
\sigma^s  =  \mathrm{a}(\alpha)(2 \mu + 3 \lambda) \varepsilon^s,
\end{eqnarray}
so that
\begin{eqnarray}
\sigma^d:\sigma^d  =  \mathrm{a}(\alpha)2 \mu \sigma^d :\varepsilon^d, &
\sigma^s:\sigma^s  =  \mathrm{a}(\alpha)(2 \mu + 3 \lambda) \sigma^s: \varepsilon^s.
\end{eqnarray}

These quantities can also expressed in terms of Young's modulus and Poisson's ratio according to the relations
\begin{eqnarray}\label{lameparameters}
\displaystyle \mu = \frac{E}{2(1+\nu)}, & \displaystyle \lambda = \frac{E\nu}{(1+\nu)(1-2\nu)}.
\end{eqnarray}

The conditions for damage can thus be written, up to normalizing by a factor $\frac{1}{E}$, so as to compare with the models of Section \ref{GSM}. 
\begin{equation}\label{damacond}
\frac{1}{E}\left(\sigma^d:\sigma^d - \frac{2}{3} \kappa \sigma^s:\sigma^s\right)> 0.
\end{equation}

We finally arrive at the form of the function $H(\varepsilon,\alpha)$ as the variation of the damage condition \eqref{damacond} with respect to $\alpha$
\begin{equation}
H(\varepsilon,\alpha) = \frac{\partial}{\partial \alpha}\left( \frac{1}{E}\left(\sigma^d:\sigma^d - \frac{2}{3} \kappa \sigma^s:\sigma^s\right)\right).
\end{equation}

In the case of two dimensions, the characteristic equation takes the following form
\begin{equation}
\mathrm{det}\left(\sigma - \lambda I\right) = \lambda^2 - \mathrm{tr}(\sigma)+ \mathrm{det}(\sigma), 
\end{equation}
where the eigenvalues of $\sigma$ are defined by
\begin{eqnarray}
\lambda_{1}=\frac{\mathrm{tr}(\sigma)+ \sqrt{\mathrm{tr}(\sigma)^2-4\mathrm{det}(\sigma)}}{2}, & \displaystyle \lambda_{2}=\frac{\mathrm{tr}(\sigma)- \sqrt{\mathrm{tr}(\sigma)^2-4\mathrm{det}(\sigma)}}{2}.
\end{eqnarray}

Then, following the Mohr's circle, we can assert that the largest shear component is obtained when
\begin{equation}
\frac{\lambda_1 - \lambda_2}{2} = \frac{1}{2}\sqrt{\mathrm{tr}(\sigma)^2 - 4 \mathrm{det}(\sigma)}.
\end{equation}

Considering $m= \frac{1}{2} \mathrm{tr}(\sigma)$, we have
\begin{eqnarray}
\vert m \vert = \sqrt{m^2} = \sqrt{\frac{1}{2} \sigma^s : \sigma^s} , &
\frac{1}{2}\sqrt{\mathrm{tr}(\sigma)^2 - 4 \mathrm{det}(\sigma)} = \sqrt{\frac{1}{2} \sum_{ij}(\sigma-mI)^2_{ij}} = \sqrt{\frac{1}{2} \sigma^d : \sigma^d},
\end{eqnarray} 
so, in two-dimensional case, damage takes place when
\begin{equation}
\sigma^d:\sigma^d - \kappa \sigma^s:\sigma^s > 0,
\end{equation}
and the function $H(\varepsilon,\alpha)$ takes the following form
\begin{equation}
H(\varepsilon,\alpha) = \frac{\partial}{\partial \alpha}\left( \frac{1}{E}\left(\sigma^d:\sigma^d - \kappa \sigma^s:\sigma^s\right)\right).
\end{equation}

We define the associated evolution by the following conditions:
\begin{enumerate} 
	\item The stress tensor $\sigma(x,t)=\sigma(u(x,t),\alpha(x,t))=\mathrm{a}(\alpha(x,t))A_0\varepsilon(u(x,t))$ satisfies the equilibrium equations
	\begin{equation}
	\begin{array}{r c l c l}
	\mathrm{div}(\sigma(x,t)) +f(x,t)  &= &0 & \mbox{ in } &\Omega, \\
	\sigma(x,t) \cdot n &= &F(x,t) &\mbox{ on } &\partial \Omega_F, \\
	u(x,t) &= &U(x,t) &\mbox{ on } &\partial \Omega_U.
	\end{array}
	\end{equation}
	
	\item The damage field $\alpha(x,t)$ satisfies the nonlocal damage criterion in $\Omega$
	\begin{equation}\label{SCdamcrit}
	\frac{(\mathrm{a}^2(\alpha(x,t)))'}{2E}\left( (A_0\varepsilon)^d:(A_0\varepsilon)^d - c (A_0\varepsilon)^s:(A_0\varepsilon)^s \right)+ w'(\alpha(x,t))-w_1\ell^2 \Delta \alpha(x,t) \geq 0,
	\end{equation}
	and the nonlocal consistency condition in $\Omega$
	\begin{equation}
	\left(\frac{(\mathrm{a}^2(\alpha(x,t)))'}{2E}\left( (A_0\varepsilon)^d:(A_0\varepsilon)^d - c (A_0\varepsilon)^s:(A_0\varepsilon)^s \right)+ w'(\alpha(x,t))-w_1\ell^2 \Delta \alpha(x,t) \right) \dot{\alpha}(x,t) = 0, 
	\end{equation}
	where $c=\kappa$ in the two-dimensional case and $c=\frac{2}{3}\kappa$ in the three-dimensional case.	
\end{enumerate}

\section{Numerical results}\label{NumSec}

\subsection{Discretization and numerical implementation}

The three models have been numerically implemented in two dimension, following \cite{marigo2016overview}. Our numerical scheme uses an alternate minimization algorithm, which consists in solving a series of subproblems at each time step, to determine $u$ when $\alpha$ is fixed, then to determine $\alpha$ at fixed $u$, until convergence.

The evolution is discretized in time. Given the displacement and the damage field $(u_{i-1},\alpha_{i-1})$ at time step $t_{i-1}$, the displacement $u_i$ at time $t_i$ is first obtained by solving the variational problem: Find $u$ such that
\begin{eqnarray}\label{elasteq} 
\forall  v \in \mathcal{C}_{t_i}, & \displaystyle \int_{\Omega}\sigma(u,\alpha_{i-1}):\varepsilon(v) dx = \int_{\Omega} f_i \cdot v dx + \int_{\partial \Omega_F} F_i \cdot v dS,
\end{eqnarray}
Subsequently, the field $\alpha_i$ is determined as a solution to the following bound-constrained minimization problem
\begin{equation}\label{moddamageeq}
\inf\left\{ \overline{\mathcal{P}}(u_i,\alpha) \mbox{ : } \alpha \in \mathcal{D}_t(\alpha_{i-1})\right\},
\end{equation}
where the functional $\overline{\mathcal{P}}(u,\alpha)$ is defined for the three different models by
\\{\it Gradient damage model:}
\begin{equation}\label{moddamageeqmod1}
\overline{\mathcal{P}}(u,\alpha)=\int_{\Omega} \frac{1}{2}\mathrm{a}(\alpha)A_0\varepsilon(u):\varepsilon(u)+w(\alpha)+\frac{1}{2}w_1\ell^2\vert \nabla \alpha \vert^2,
\end{equation}
{\it Gradient damage model for shear fracture:}
\begin{equation}\label{moddamageeqmod2}
\overline{\mathcal{P}}(u,\alpha)=\int_{\Omega} \left(\lambda + \frac{2 \mu}{n} \right) \frac{\mathrm{tr}(\varepsilon(u))^2}{2}+\mathrm{a}(\alpha)\mu \varepsilon^d:\varepsilon^d + w(\alpha) + \frac{1}{2} w_1\ell^2 \vert \nabla\alpha \vert^2,
\end{equation}
{\it Shear-compression damage model:}
\begin{equation}\label{moddamageeqmod3}
\overline{\mathcal{P}}(u,\alpha)=\int_{\Omega} \frac{\mathrm{a}^2(\alpha)}{2E} \left( \left(A_0\varepsilon(u)\right)^d:\left(A_0\varepsilon(u)\right)^d-\kappa \left(A_0\varepsilon(u)\right)^s:\left(A_0\varepsilon(u)\right)^s \right)+w(\alpha)+w_1\ell^2 \vert \nabla \alpha\vert^2.
\end{equation}

The unilateral constraint $\alpha(x)\geq \alpha_{i-1}$ is the time-discrete version of the irreversibility of damage. The solution strategy is summarized in Algorithm \ref{alg:damage}. The unilateral constraint included in the damage problem requires the use of variational inequalities solvers, here we use the open-source library PETSc \cite{balay2004petsc}, where this capability is available. The problem is discretized in space with standard triangular finite elements with piecewise linear approximation for $u$ and $\alpha$. The mesh size is selected so that localization band contains at least 3 elements. Finite element implementations based on the open-source FEniCS library \cite{langtangen2016solving, logg2012automated} have been used.
\begin{algorithm}[H]
\begin{algorithmic}[1]
	\RETURN Solution of time step $t_i$.
	\STATE Given $(u_{i-1},\alpha_{i-1})$, the sate at the previous loading step.
	\STATE Set $(u^{(0)},\alpha^{(0)}):=(u_{i-1},\alpha_{i-1})$ and error$^{(0)} = 1.0$ 
	\WHILE {error$^{(p)}>$ tolerance}
	\STATE Solve $u^{(p)}$ from (\ref{elasteq}) with $\alpha^{(p-1)}$.
	\STATE Find $\displaystyle\alpha^{(p)}:=  \argmin_{\alpha \in \mathcal{D}(\alpha_{i-1})} \overline{\mathcal{P}}(u^{(p)},\alpha)$.
	\STATE error$^{(p)} = \Vert\alpha^{(p-1)} -\alpha^{(p)}\Vert_{\infty}$.
	\ENDWHILE
	\STATE  Set $(\mathbf{u}_{i}, \alpha_{i}) = (u^{p}, \alpha^{p}).$
\end{algorithmic}
\caption{Numerical algorithm to solve the damage problem}\label{alg:damage}
\end{algorithm}

\subsection{Influence of the cavity in the damage model}

\begin{figure}[h!]
	\centering
	\includegraphics[scale=1.0]{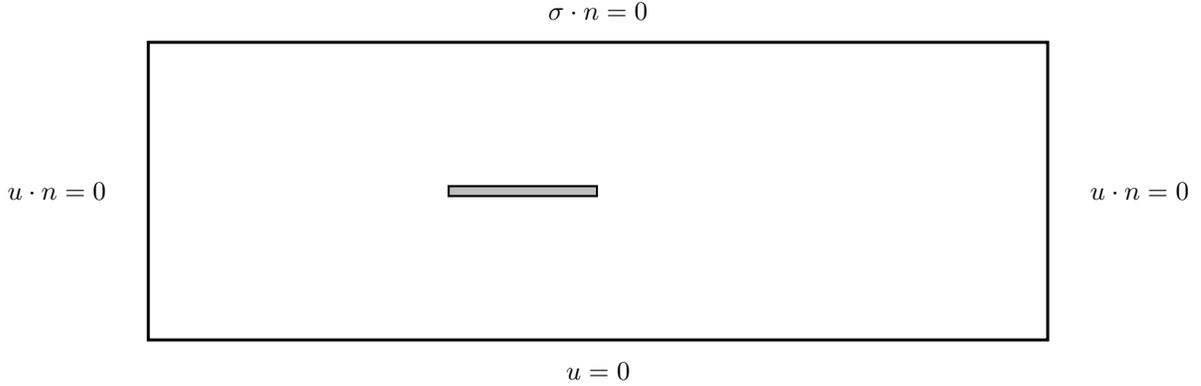}
	\caption{Geometry and boundary condition for the cavity problem.}
	\label{geometry}
\end{figure}

For modeling the block caving process, we consider a domain that simulates a rock mass in two dimensions $\Omega_0\subseteq \mathbb{R}^2$ with $\partial \Omega_0 = \Gamma_{lat} \cup \Gamma_{up} \cup \Gamma_{down}$ denoting the lateral, upper and lower bounded of the domain respectively. We consider a time depending interior domain $S(t_i)\subseteq \Omega_0$, with $\partial S(t_i)=\Gamma_{cav}(t_i)$ representing the interior cavity of this rock mass, where we assume free boundary conditions in $\Gamma_{cav}(t_i)$, that is, $\sigma \cdot n =0$ in $\Gamma_{cav}(t_i)$. Finally, the domain to consider the block caving process is $\Omega(t_i) = \Omega_0 \setminus S(t_i)$. Figure \ref{geometry} shows a sketch of the domain defined above with the boundary conditions used in these problems. The loading is given by the gravity defined by $f_t=\rho g$, with $\rho=2.7\cdot 10^3\left[ \frac{Kg}{m^3}\right]$ and $g= \left(0.0,-9.8\left[ \frac{m}{s^2}\right]\right)$.

In our tests, we consider that the rock mass is defined by $\Omega_0 = (-1500,1500)\times (-500,500)$ and the cavities are represented by $S(t)=(-500,-500+40t_i)\times(-20,20)$. All simulations presented here use the following values for material parameters
\begin{equation}
\begin{array}{c c c c c}
E =2.9 \cdot 10^{10}[Pa], & \kappa=1.0, &  \mbox{and}  & \nu = 0.3,
\end{array}
\end{equation}
where we use a quadratic damage model defined by
\begin{eqnarray}
\mathrm{a}(\alpha)=(1-\alpha)^2, &  w(\alpha)=w_1\alpha^2.
\end{eqnarray}

\subsubsection{Gradient damage model}

\begin{figure}[h!]
	\centering
	\begin{subfigure}[b]{0.45\textwidth}
		\centering
		\includegraphics[width=\textwidth]{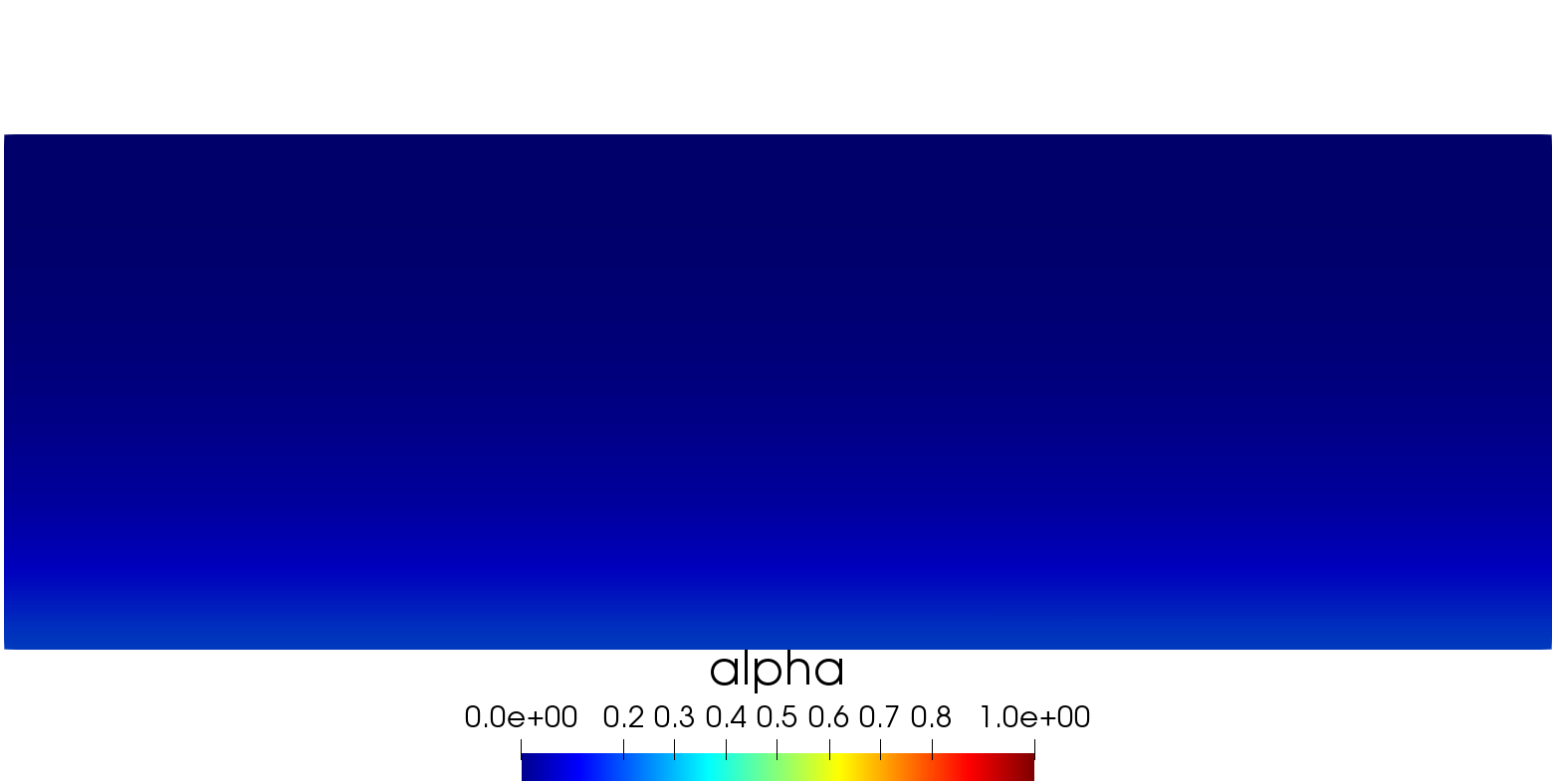}
		\caption{$t_i=0$}
	\end{subfigure}
	\begin{subfigure}[b]{0.45\textwidth}
		\centering
		\includegraphics[width=\textwidth]{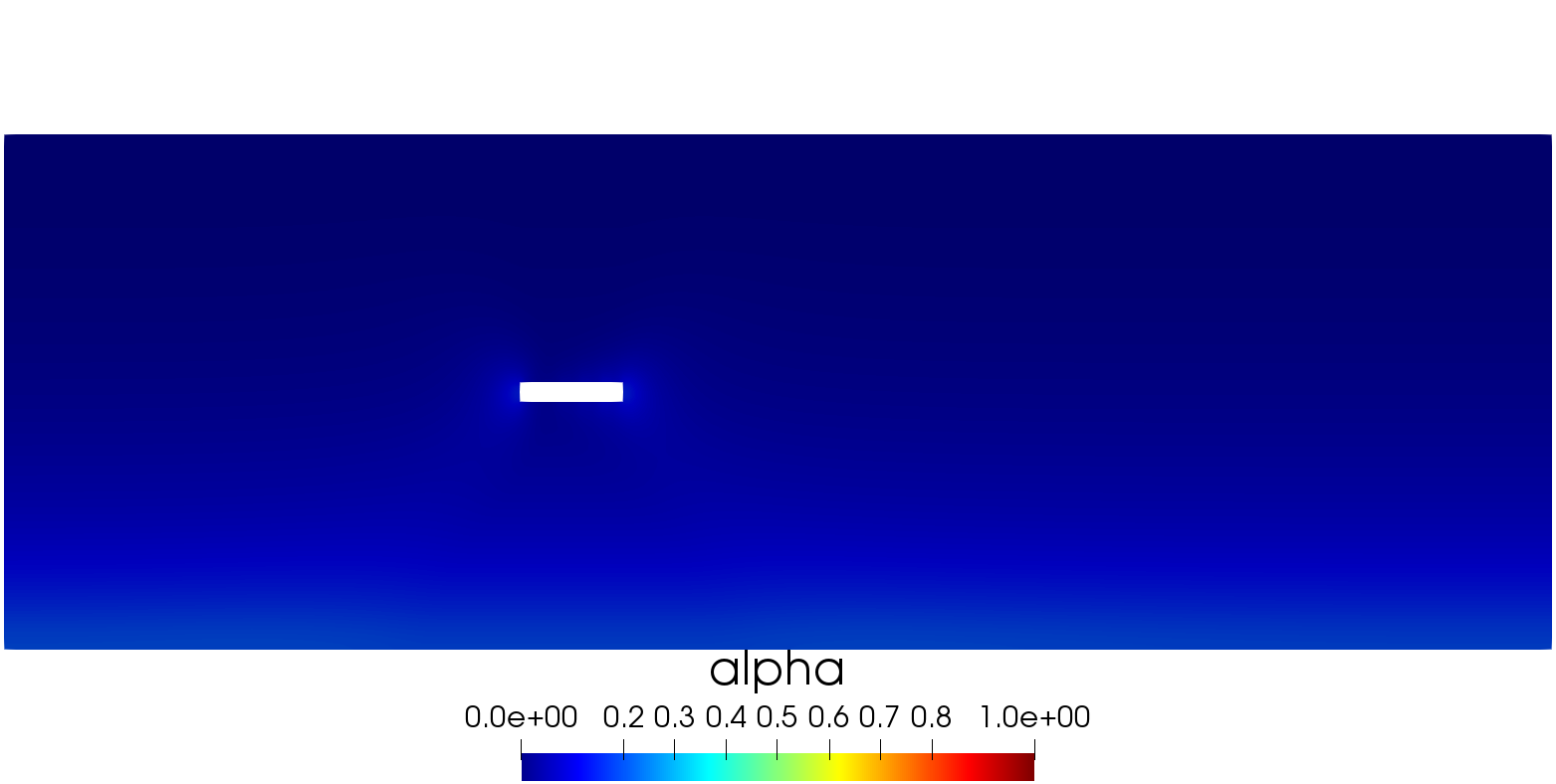}
		\caption{$t_i=5$}
	\end{subfigure}
	\begin{subfigure}[b]{0.45\textwidth}
		\centering
		\includegraphics[width=\textwidth]{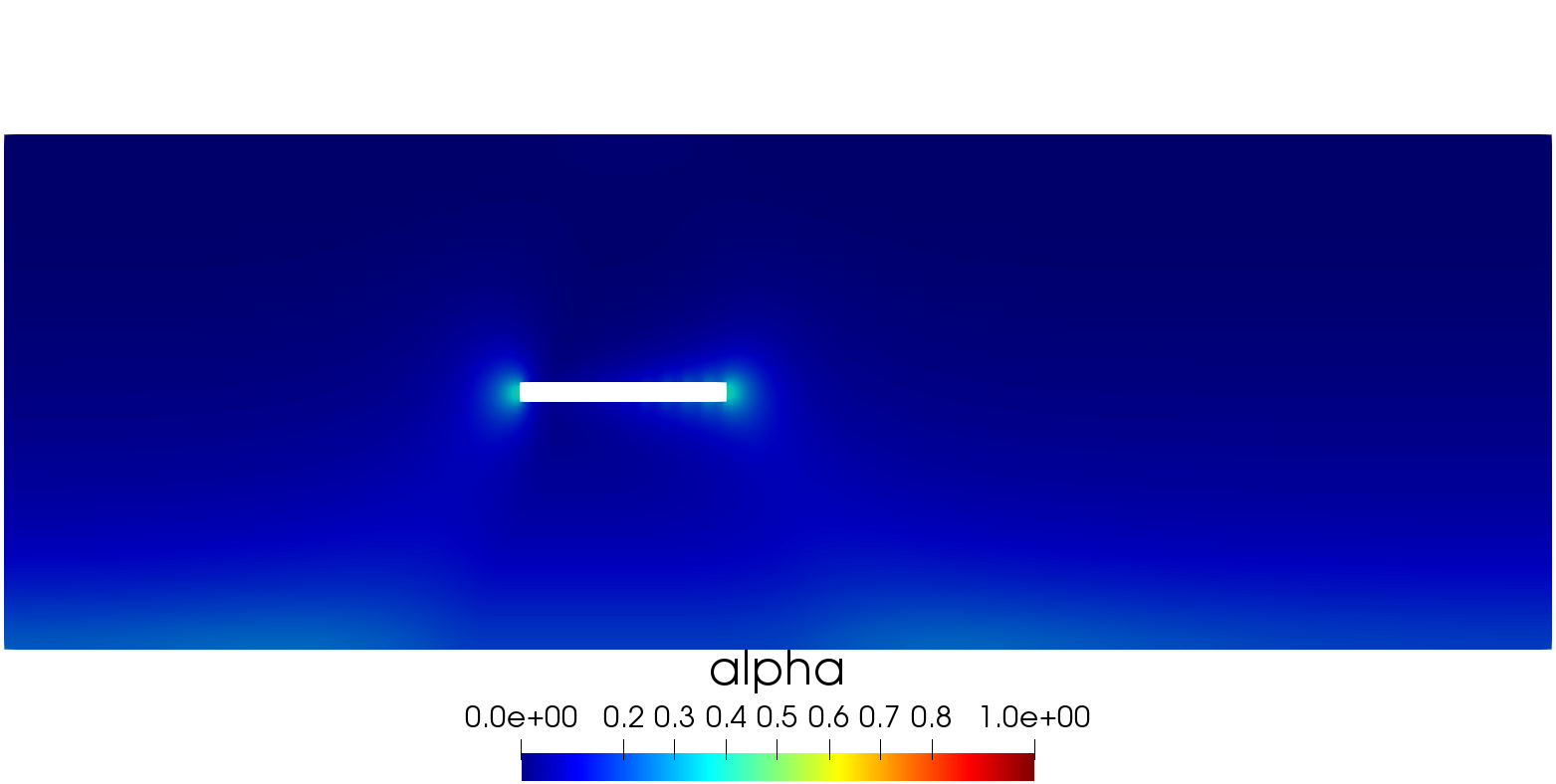}
		\caption{$t_i=10$}
	\end{subfigure}	
	\begin{subfigure}[b]{0.45\textwidth}
		\centering
		\includegraphics[width=\textwidth]{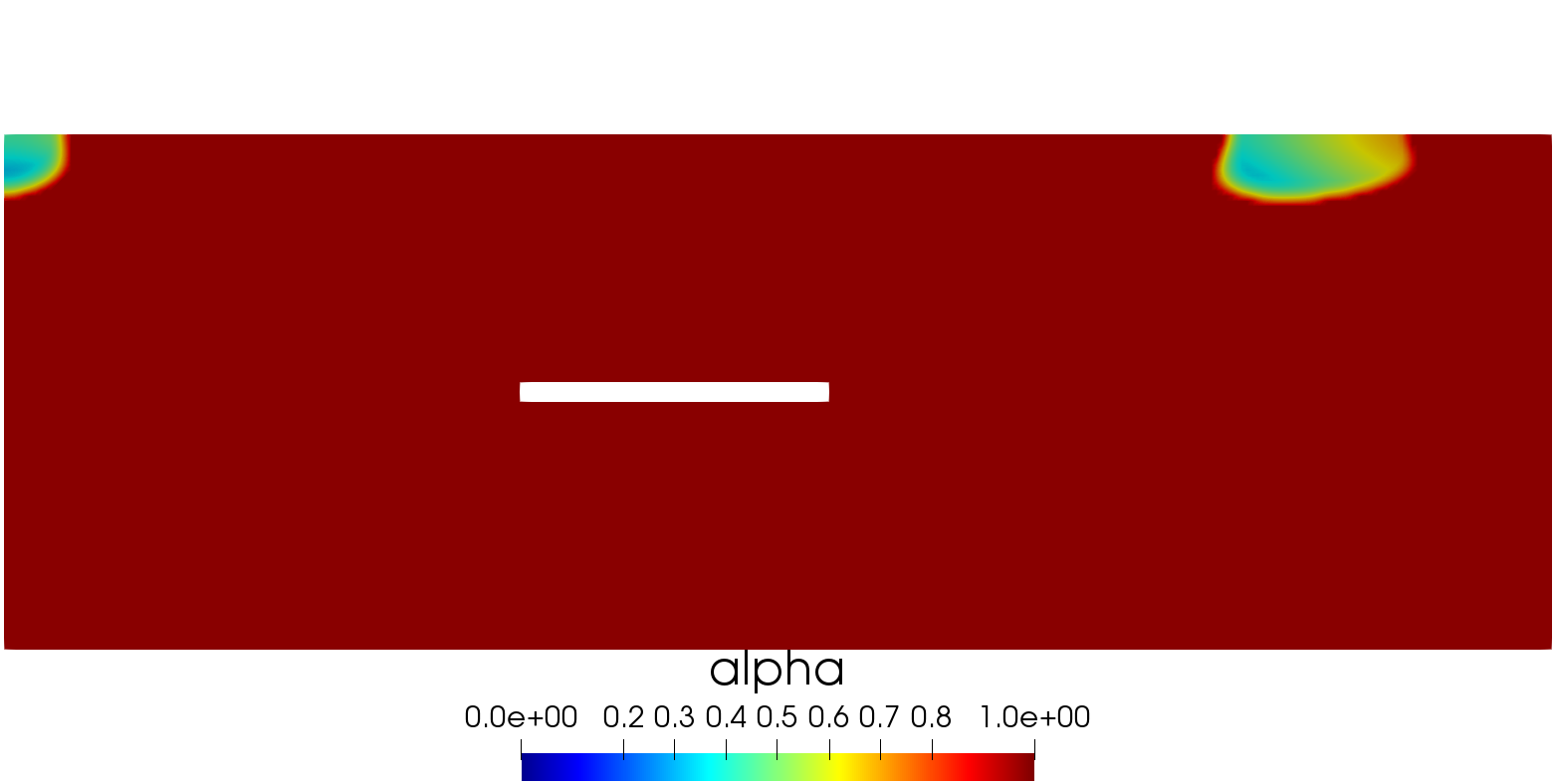}
		\caption{$t_i=15$}
	\end{subfigure}
	\caption{Damage field distribution in the rock mass for $w_1= 10^5\left[\frac{N}{m^3}\right]$.}\label{fig1}
\end{figure}

In the first test case, we consider the gradient damage model described the Section \ref{GDM}. Figures \ref{fig1} summarizes the results obtained with this model. The images display the evolution of damage when the cavity advances in time. We chose $w_1=10^5\left[\frac{N}{m^3}\right]$. We observe that the damage appears in whole the domain and it is distributed in an instant of time (the same phenomenon can be seen when considering values smaller of $w_1$). This indicates that in this example, damage is essentially triggered by compression due to the gravity forces, which is the main forcing term in the model and which is uniformly distributed in the domain. Damage does not seem to be very sensitive to the presence of the cavity.

We observe that on the first time steps, the damage begins to appear in the lower area of the domain, below the extraction cavity, and later distributes throughout the domain. Such qualitative behavior is not consistent with the expected mechanism of fracking produced by block caving.

\subsubsection{Gradient damage model for shear fracture}

\begin{figure}[h!]
	\centering
	\begin{subfigure}[b]{0.45\textwidth}
		\centering
		\includegraphics[width=\textwidth]{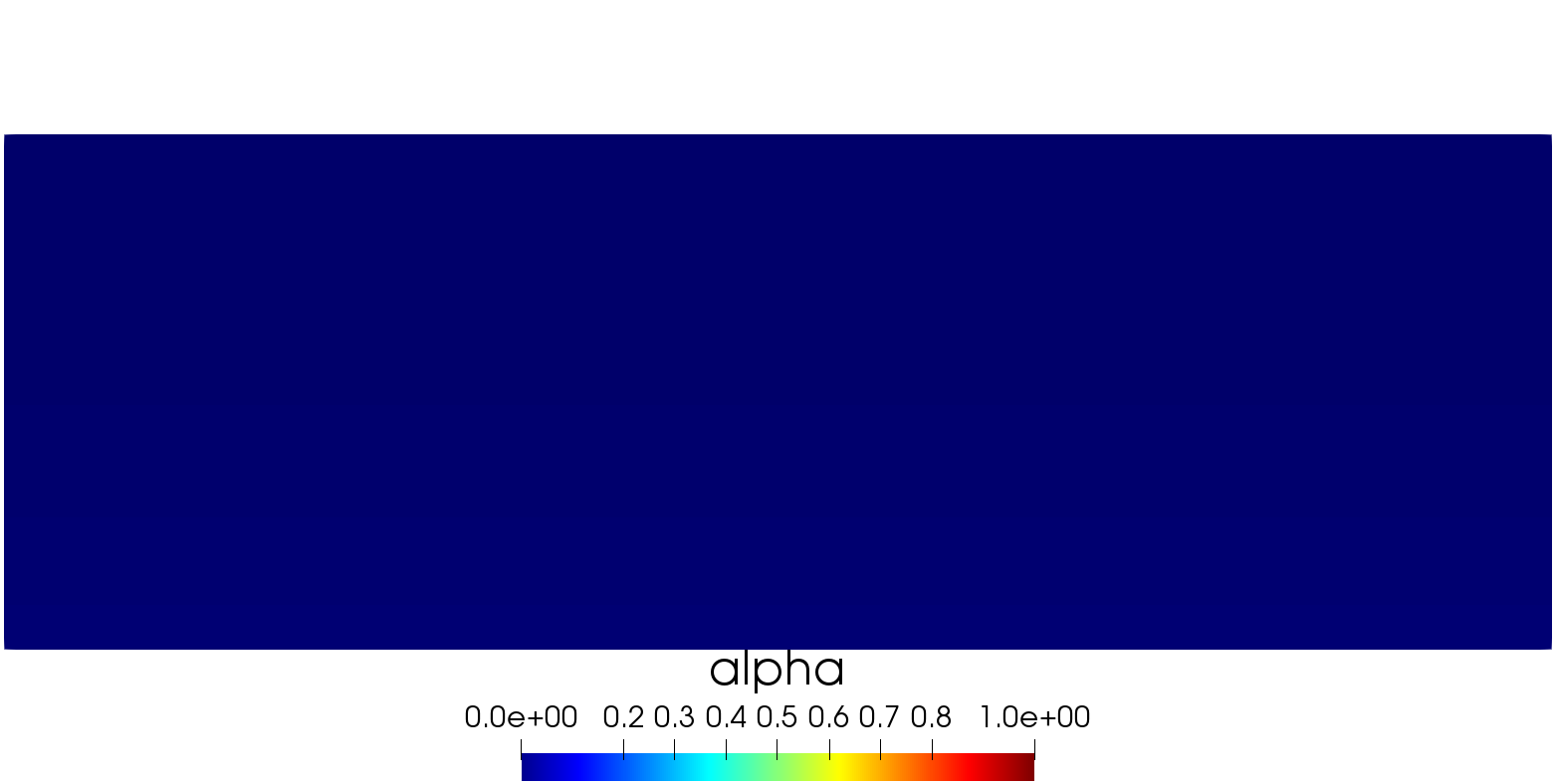}
		\caption{$t_i=0$}
	\end{subfigure}
	\begin{subfigure}[b]{0.45\textwidth}
		\centering
		\includegraphics[width=\textwidth]{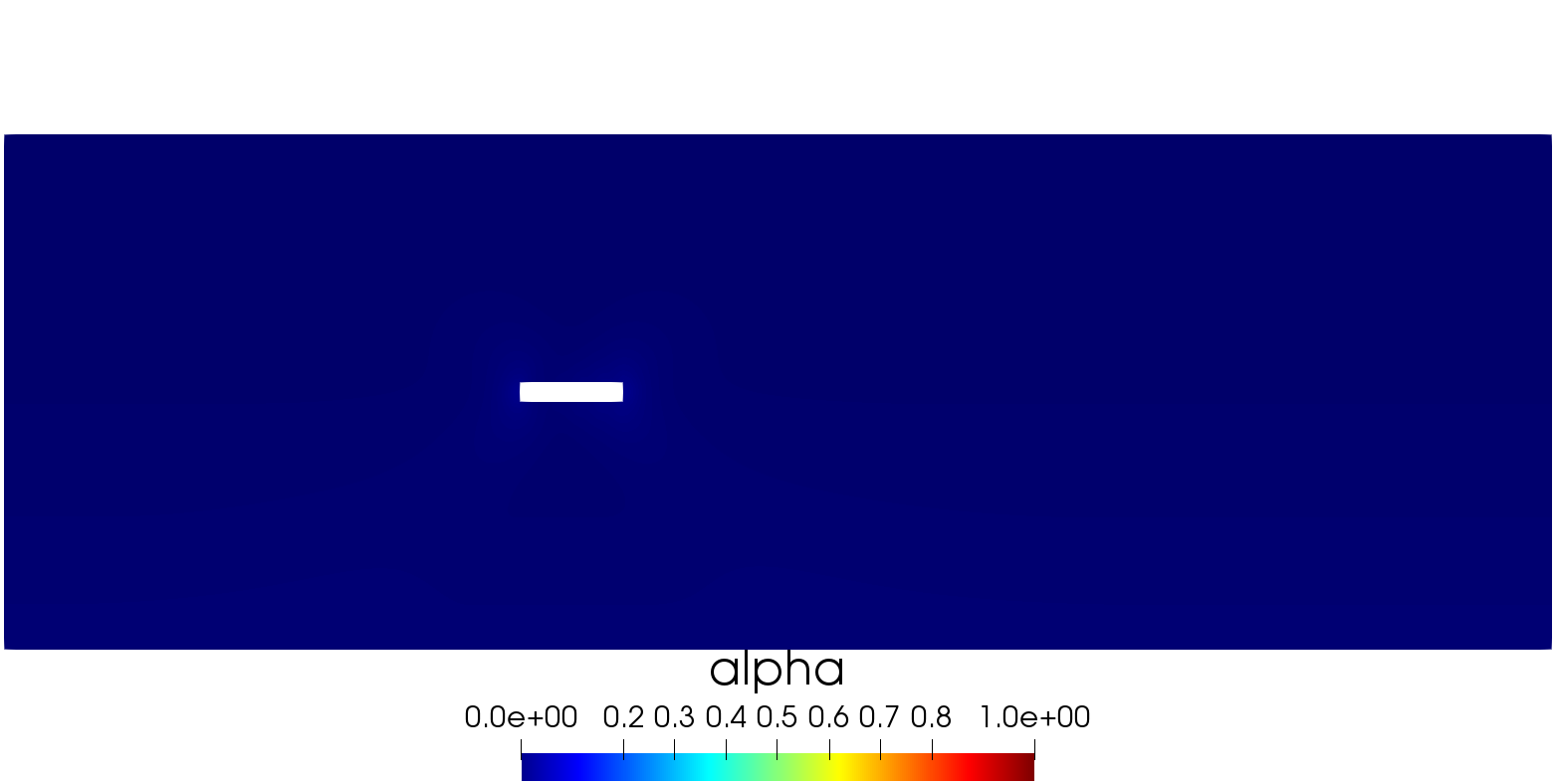}
		\caption{$t_i=5$}
	\end{subfigure}
	\begin{subfigure}[b]{0.45\textwidth}
		\centering
		\includegraphics[width=\textwidth]{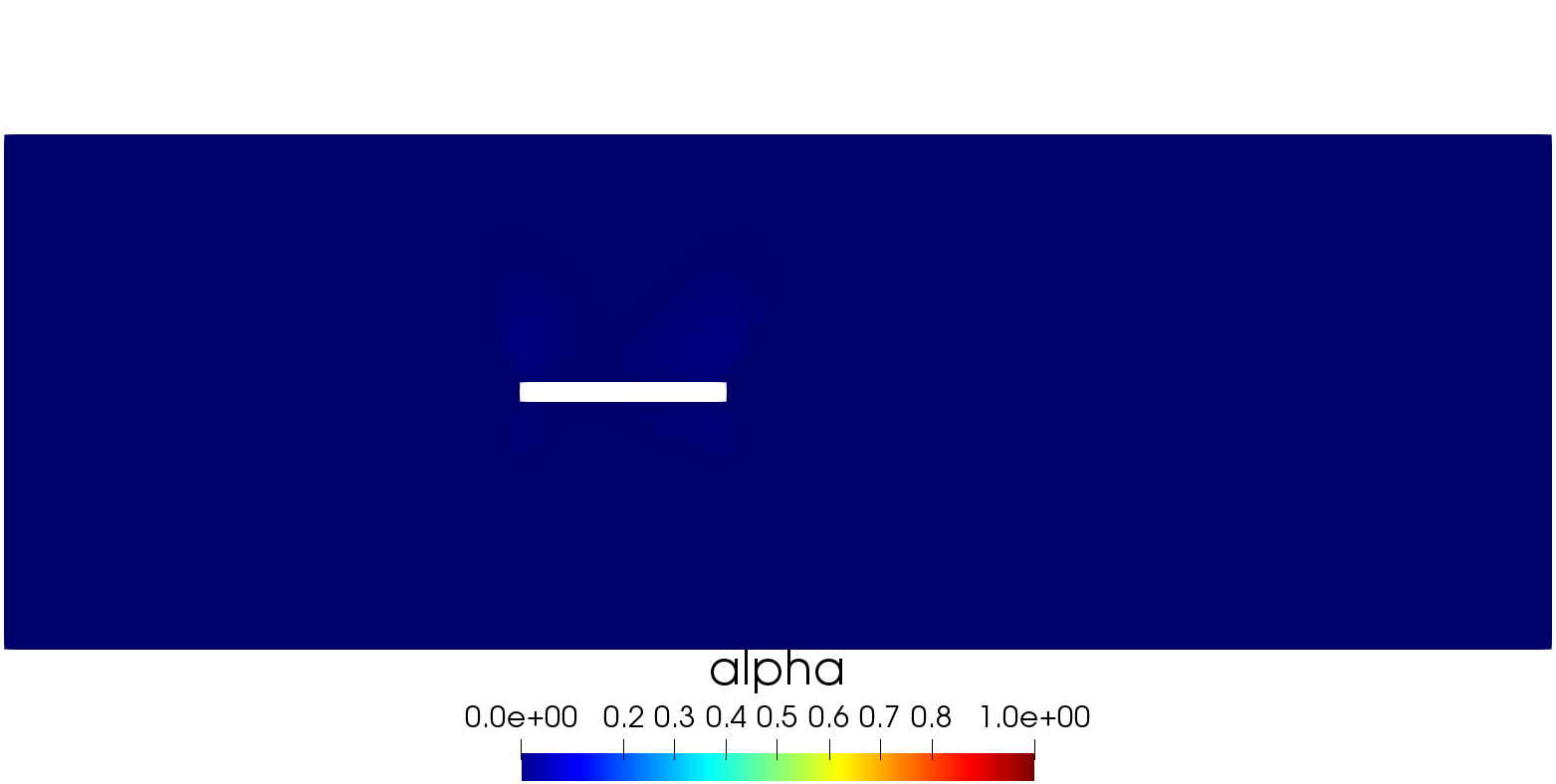}
		\caption{$t_i=10$}
	\end{subfigure}	
	\begin{subfigure}[b]{0.45\textwidth}
		\centering
		\includegraphics[width=\textwidth]{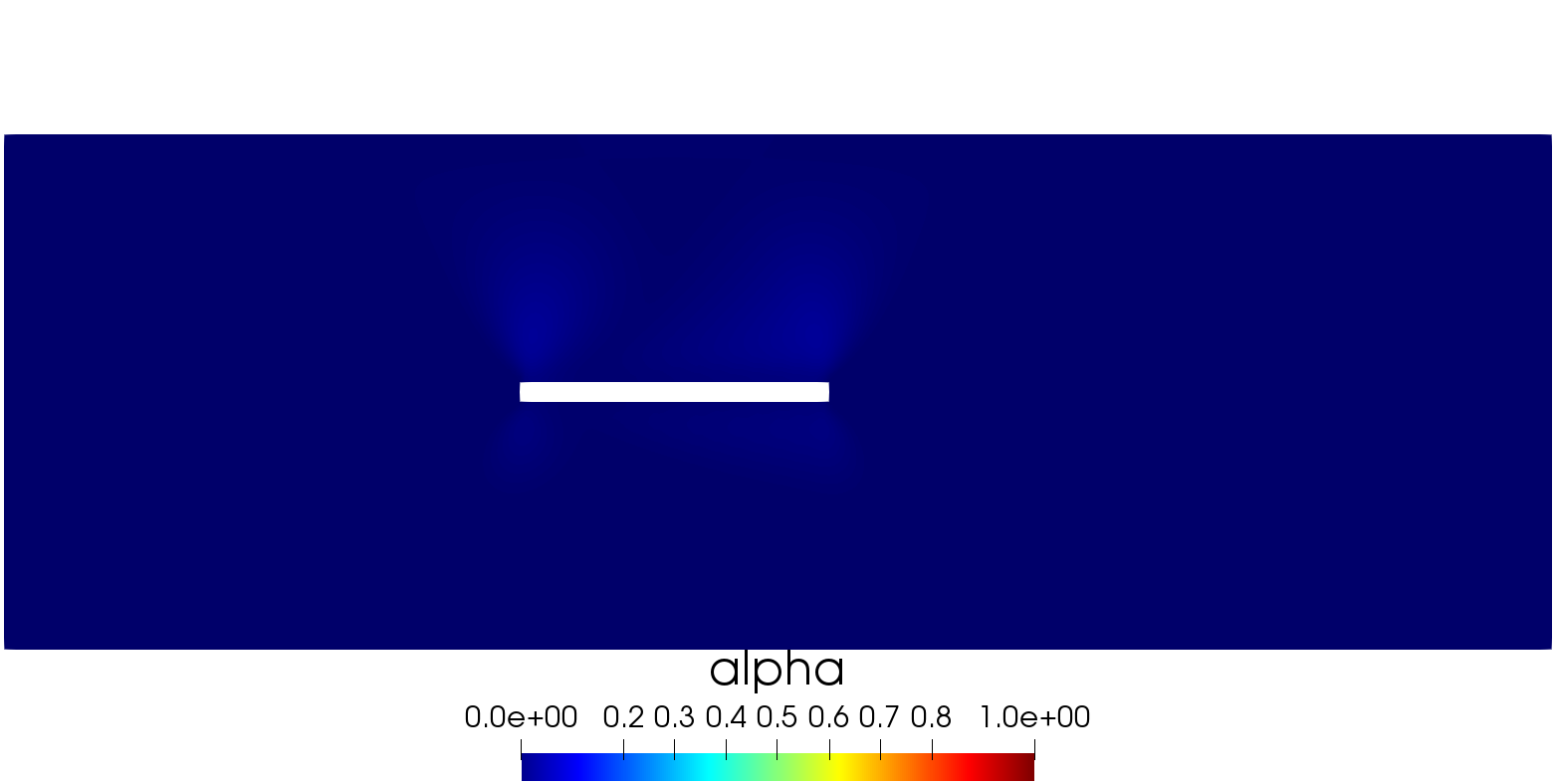}
		\caption{$t_i=15$}
	\end{subfigure}
	\caption{Damage field distribution in the rock mass for $w_1= 10^5\left[\frac{N}{m^3}\right]$}\label{fig2}
\end{figure}

In the second test case, we consider the gradient damage model for shear fracture given in Section \ref{SDM}. Figure \ref{fig2} displays the evolution of the damage when the cavity advances in time. Again, $w_1=10^5\left[\frac{N}{m^3}\right]$. Here, damage hardly appears, and takes values close to zero around the cavity.

\begin{figure}[h!]
	\centering
	\begin{subfigure}[b]{0.45\textwidth}
		\centering
		\includegraphics[width=\textwidth]{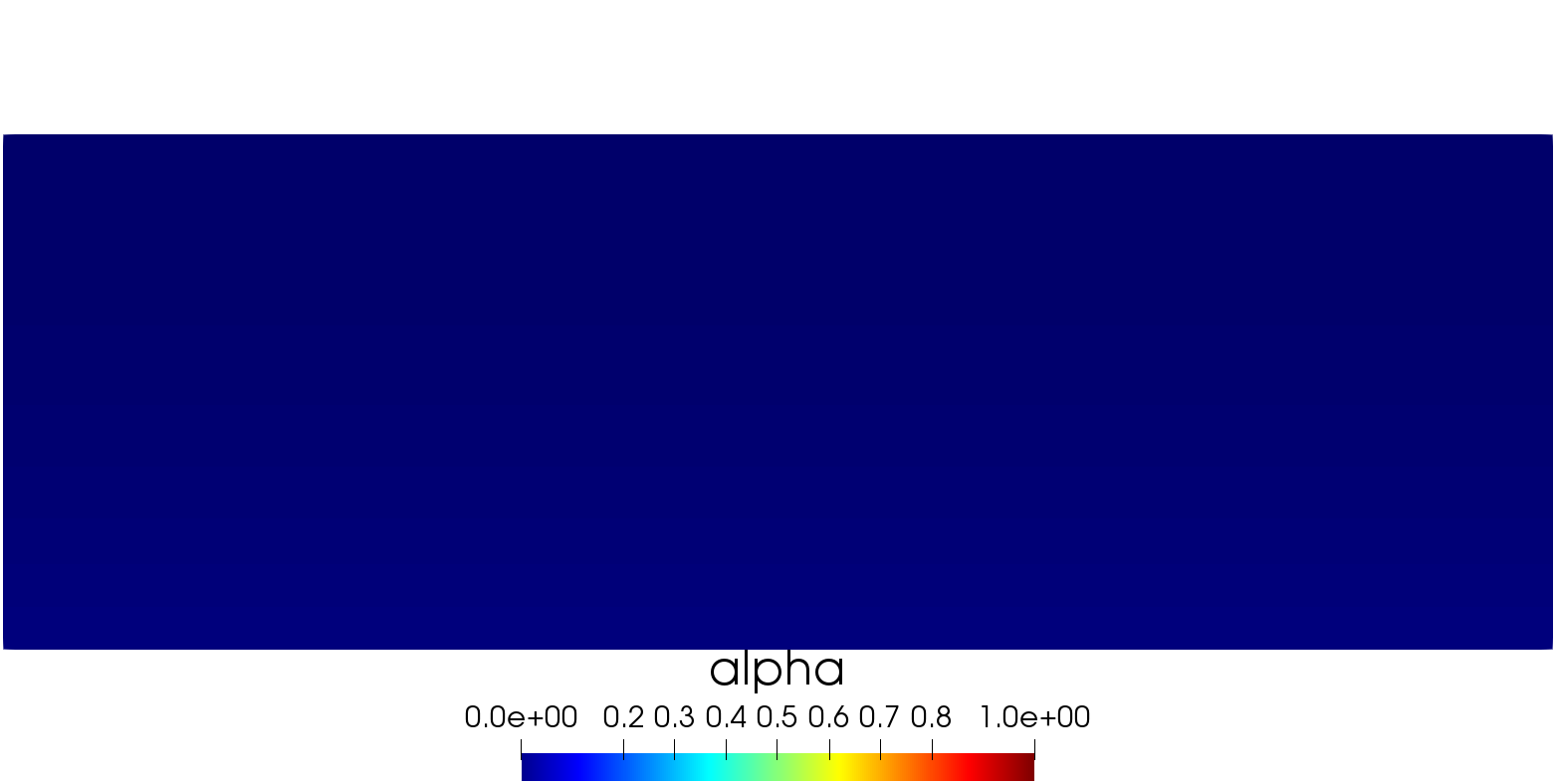}
		\caption{$t_i=0$}
	\end{subfigure}
	\begin{subfigure}[b]{0.45\textwidth}
		\centering
		\includegraphics[width=\textwidth]{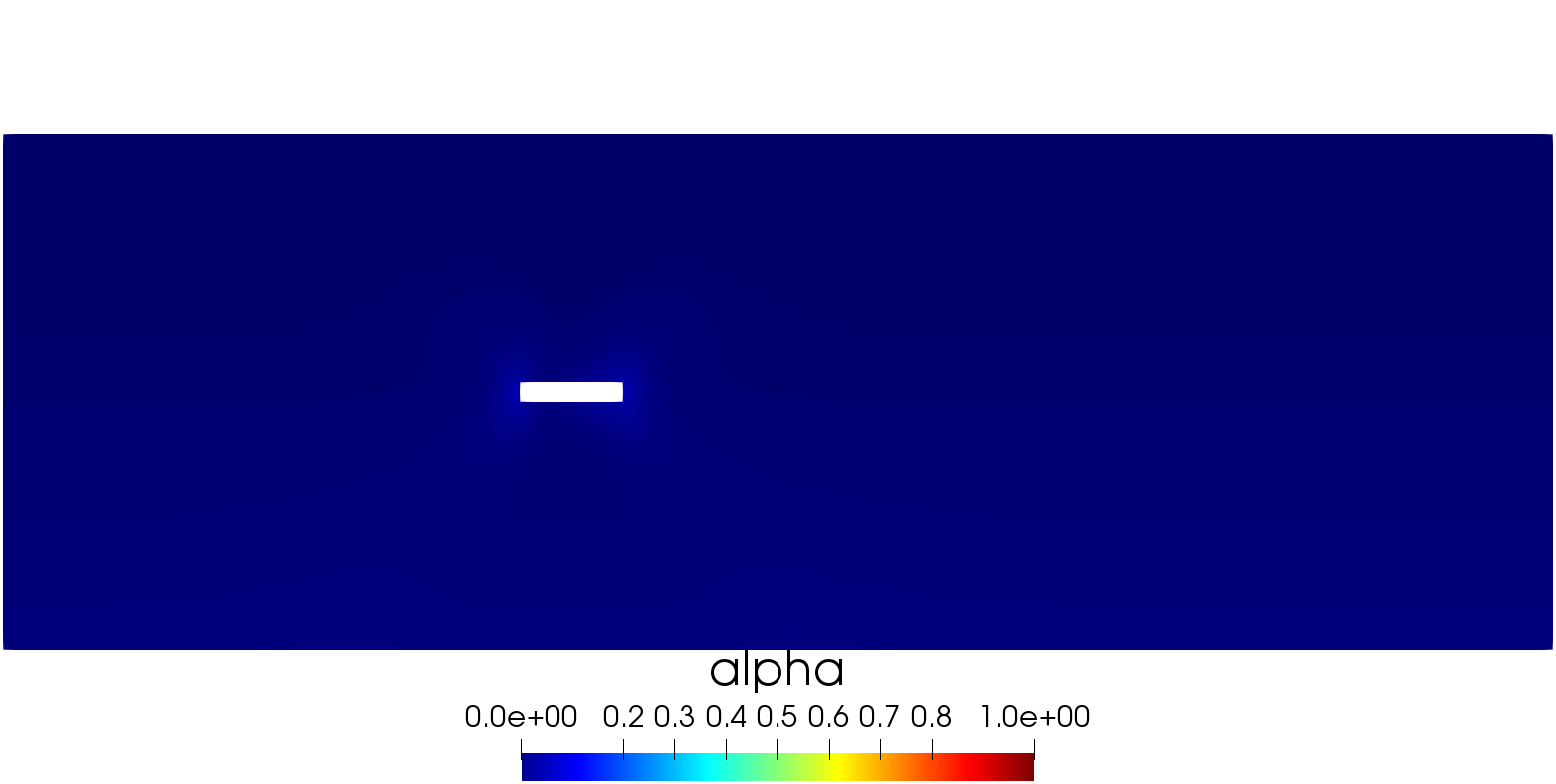}
		\caption{$t_i=5$}
	\end{subfigure}
	\begin{subfigure}[b]{0.45\textwidth}
		\centering
		\includegraphics[width=\textwidth]{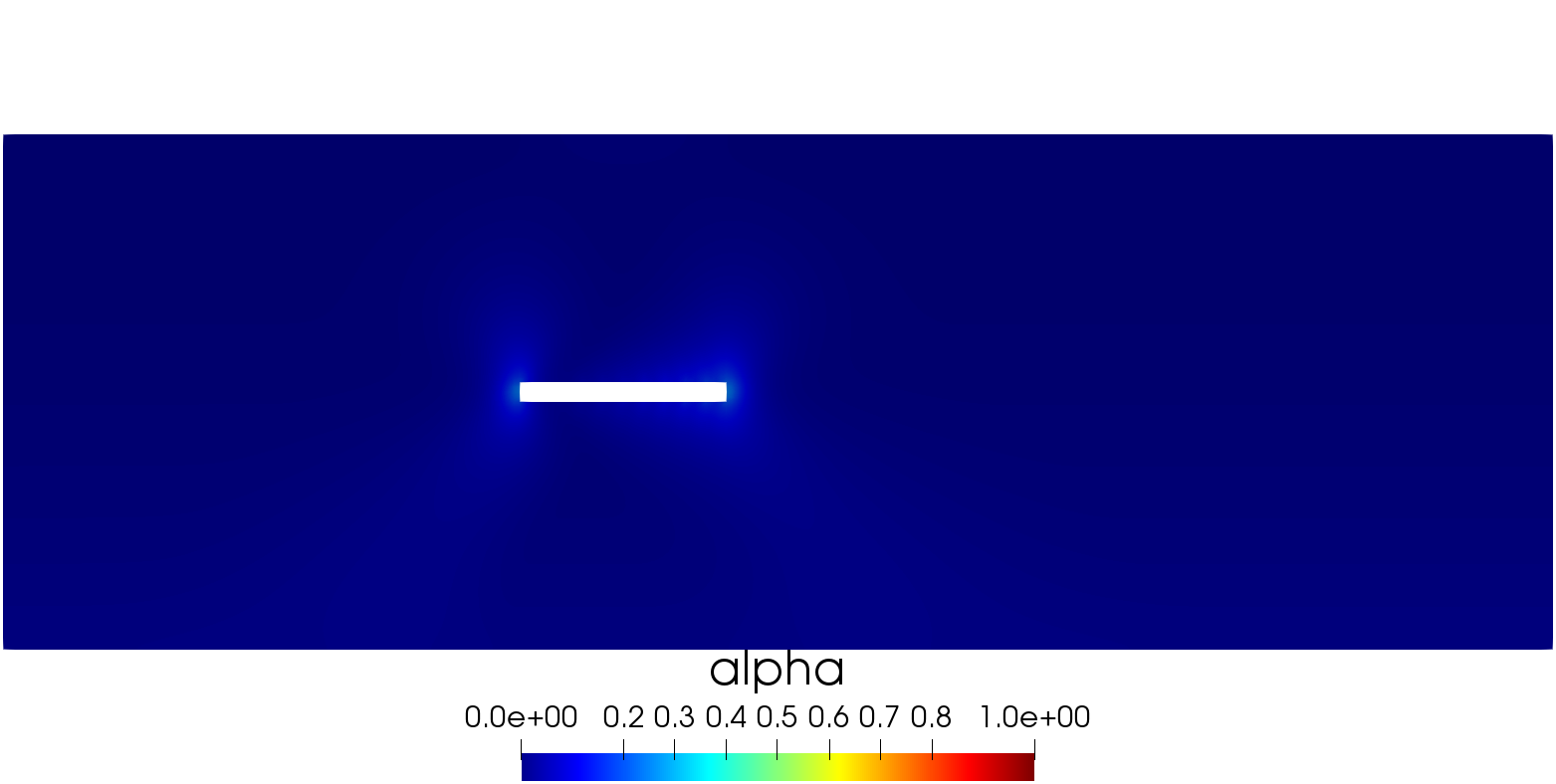}
		\caption{$t_i=10$}
	\end{subfigure}	
	\begin{subfigure}[b]{0.45\textwidth}
		\centering
		\includegraphics[width=\textwidth]{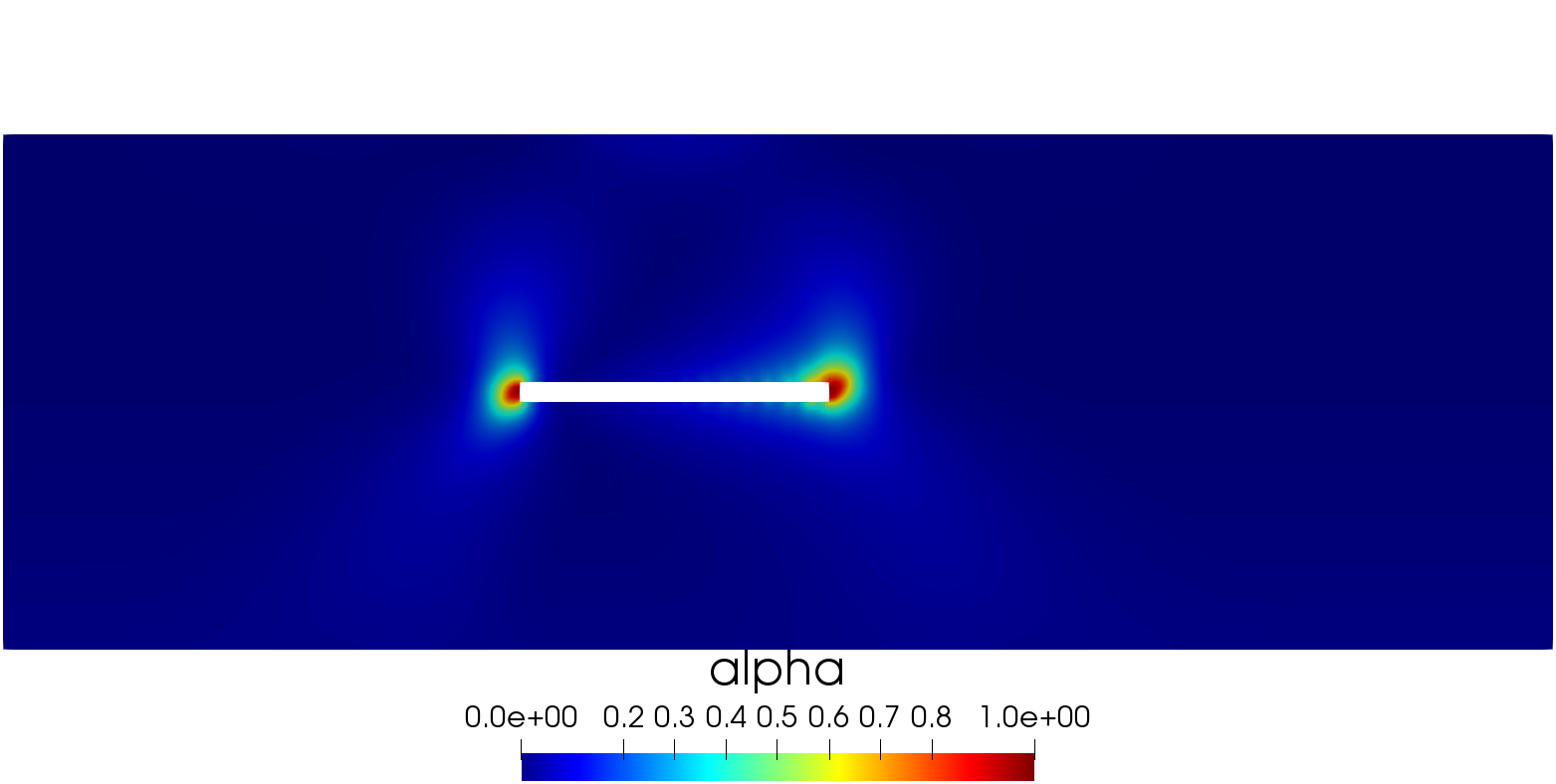}
		\caption{$t_i=15$}
	\end{subfigure}
	\caption{Damage field distribution in the rock mass for $w_1=5\cdot 10^4\left[\frac{N}{m^3}\right]$}\label{fig3}
\end{figure}

\begin{figure}[h!]
	\centering
	\begin{subfigure}[b]{0.45\textwidth}
		\centering
		\includegraphics[width=\textwidth]{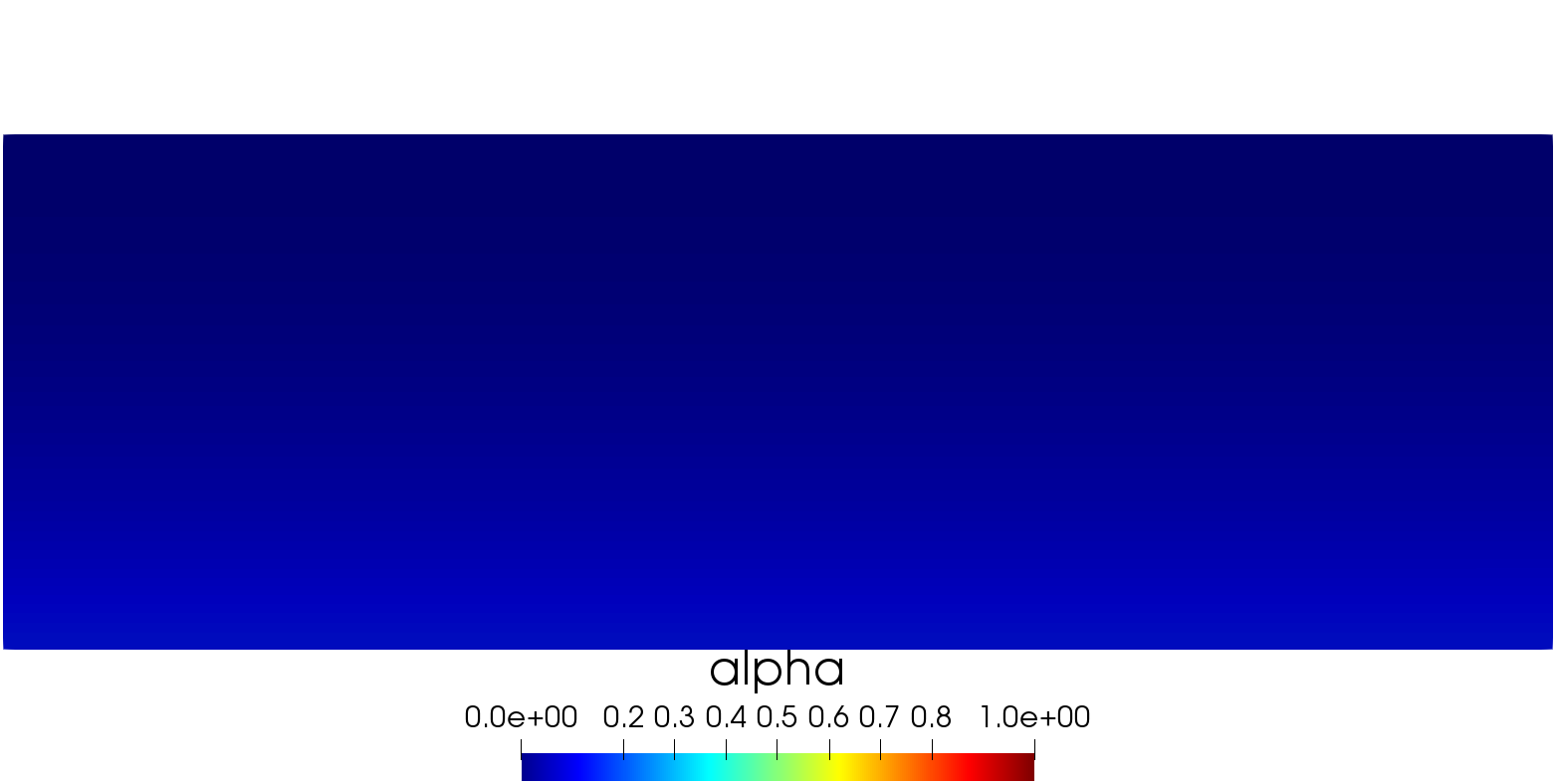}
		\caption{$t_i=0$}
	\end{subfigure}
	\begin{subfigure}[b]{0.45\textwidth}
		\centering
		\includegraphics[width=\textwidth]{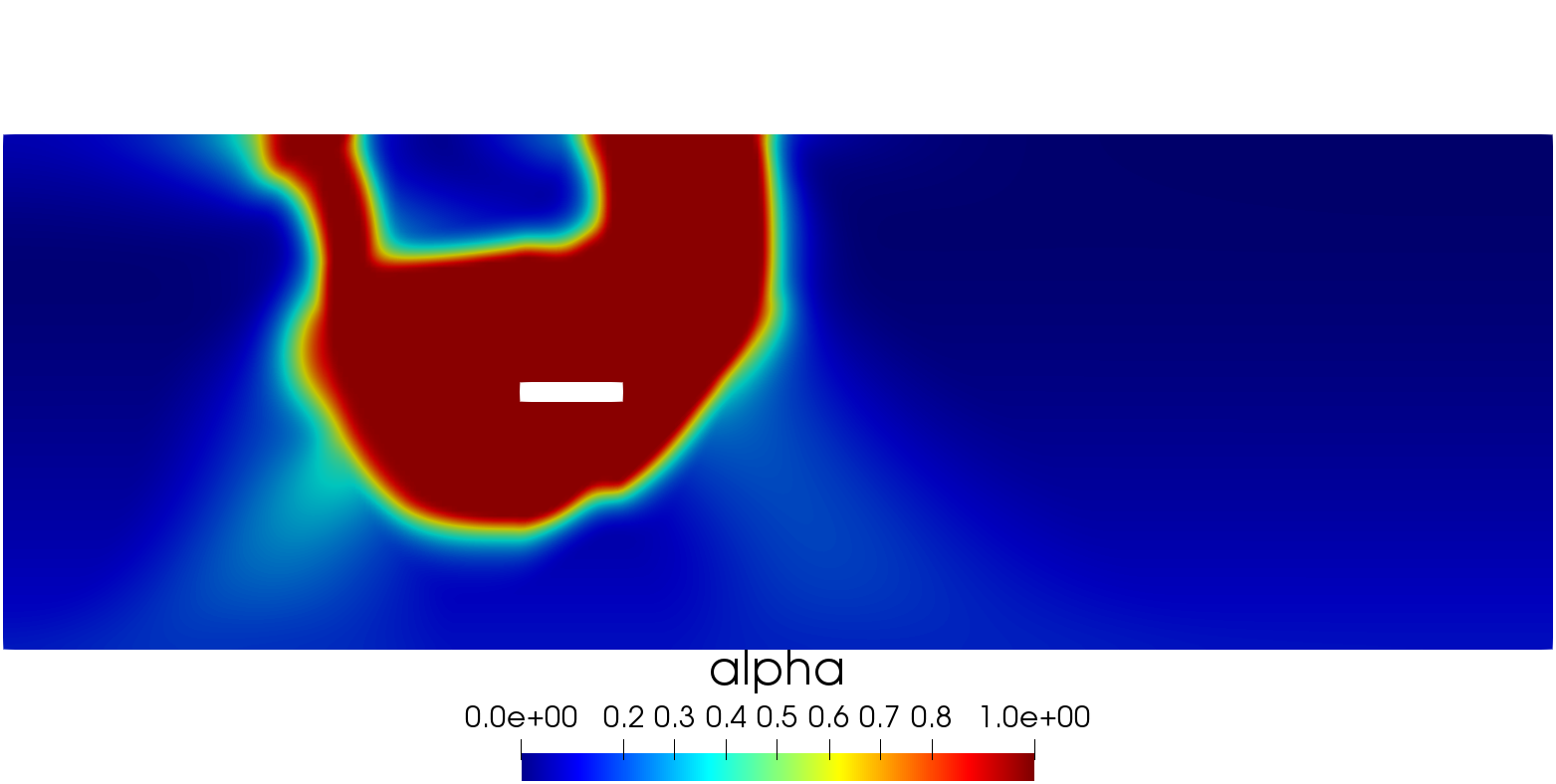}
		\caption{$t_i=5$}
	\end{subfigure}
	\begin{subfigure}[b]{0.45\textwidth}
		\centering
		\includegraphics[width=\textwidth]{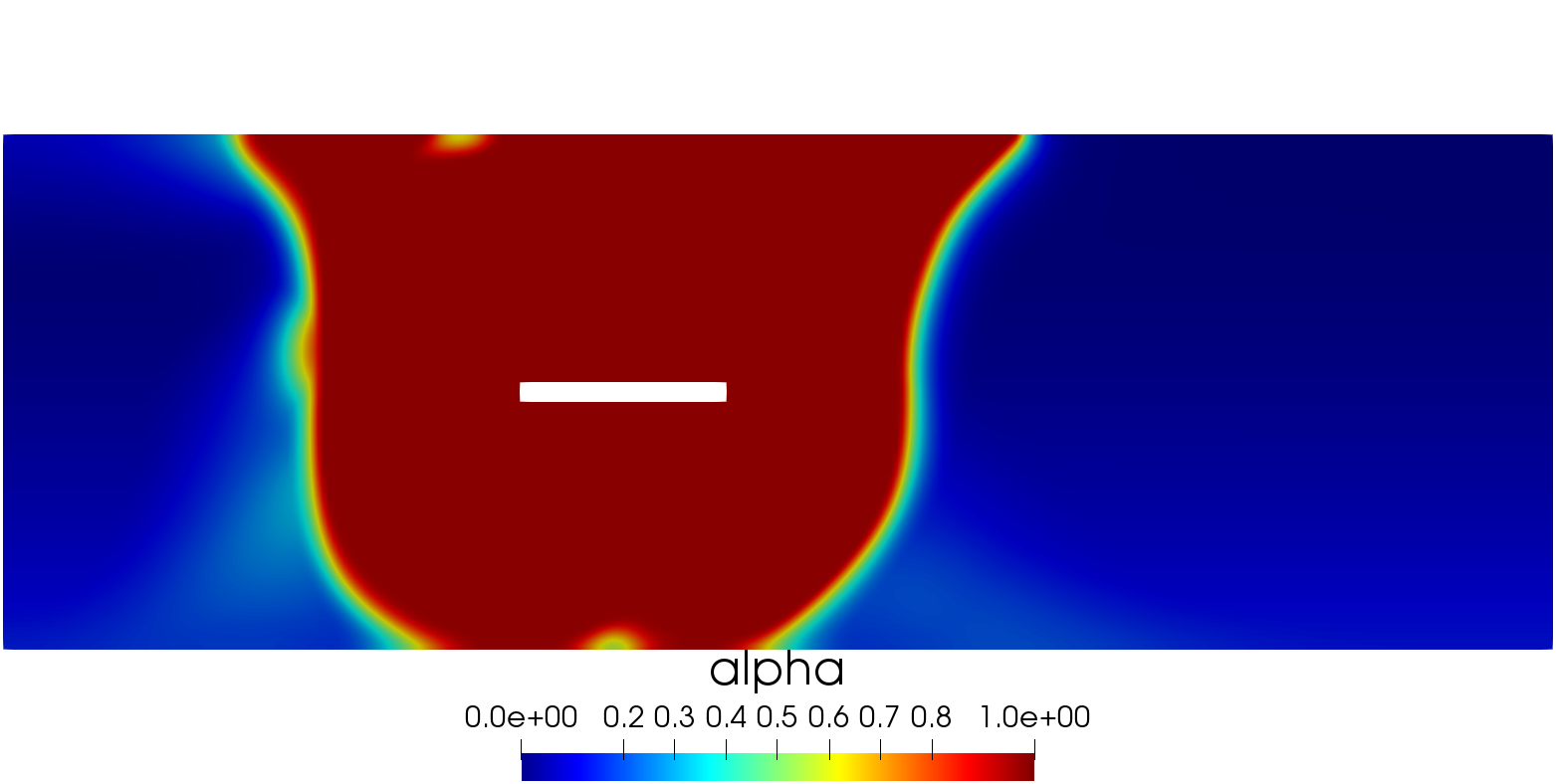}
		\caption{$t_i=10$}
	\end{subfigure}	
	\begin{subfigure}[b]{0.45\textwidth}
		\centering
		\includegraphics[width=\textwidth]{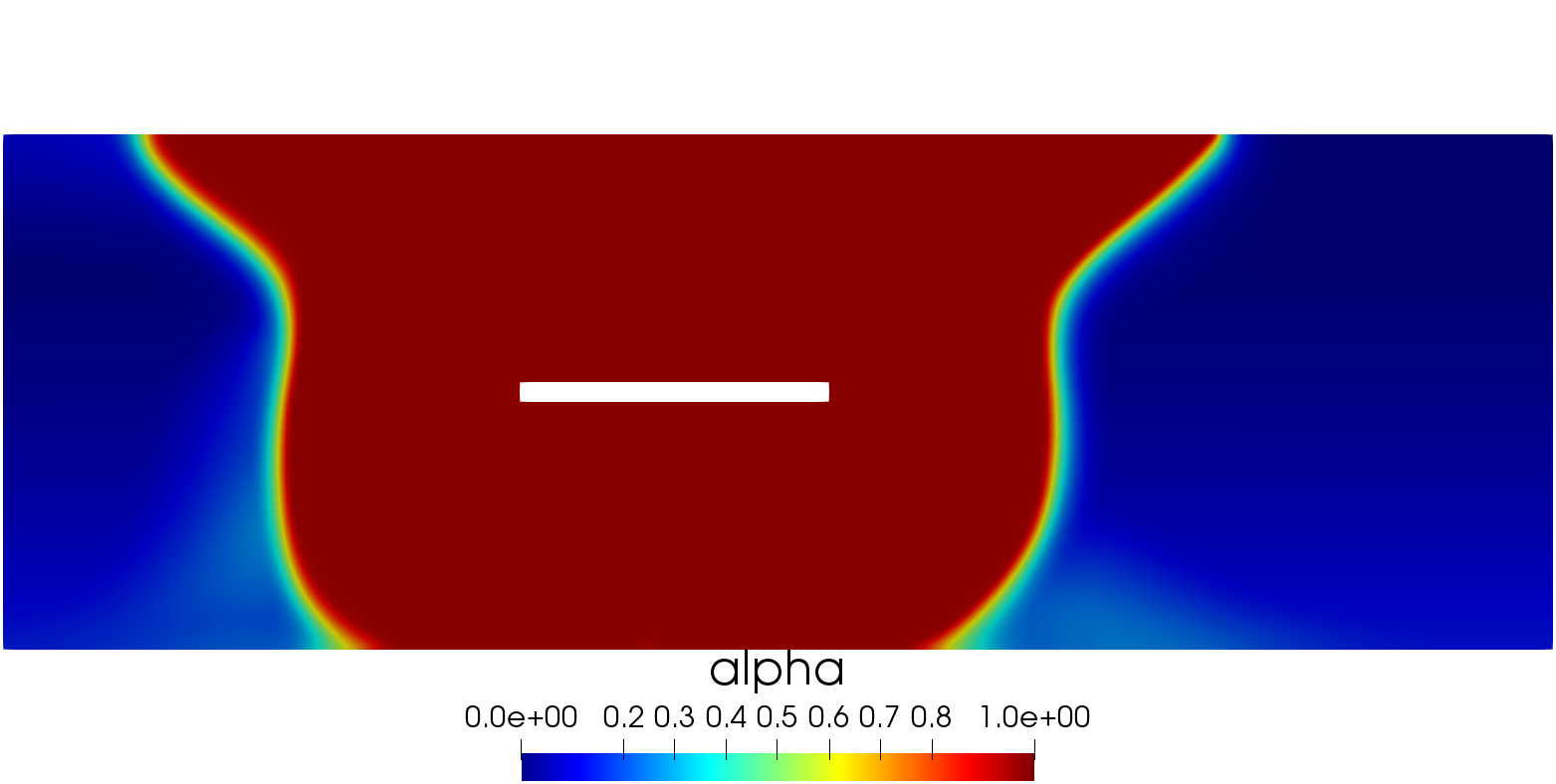}
		\caption{$t_i=15$}
	\end{subfigure}
	\caption{Damage field distribution in the rock mass for $w_1=10^4\left[\frac{N}{m^3}\right]$}\label{fig4}
\end{figure}

For smaller values of $w_1$ (see Figure \ref{fig3} and Figure \ref{fig4} for $w_1=5 \cdot 10^4\left[\frac{N}{m^3}\right]$ and $w_1= 10^4\left[\frac{N}{m^3}\right]$) damage remains distributed around the cavity. When $w_1=5 \cdot 10^4\left[\frac{N}{m^3}\right]$ begins to appear around the end corners of the cavity (where stresses are expected to blow up, due to corner singularities). On the other hand when $w_1= 10^4\left[\frac{N}{m^3}\right]$ damage begins to appear in the bottom of the rock mass. Again these results do not represent the expected effect of block caving.

\subsubsection{Shear-compression damage model}

\begin{figure}[h!]
	\centering
	\begin{subfigure}[b]{0.45\textwidth}
		\centering
		\includegraphics[width=\textwidth]{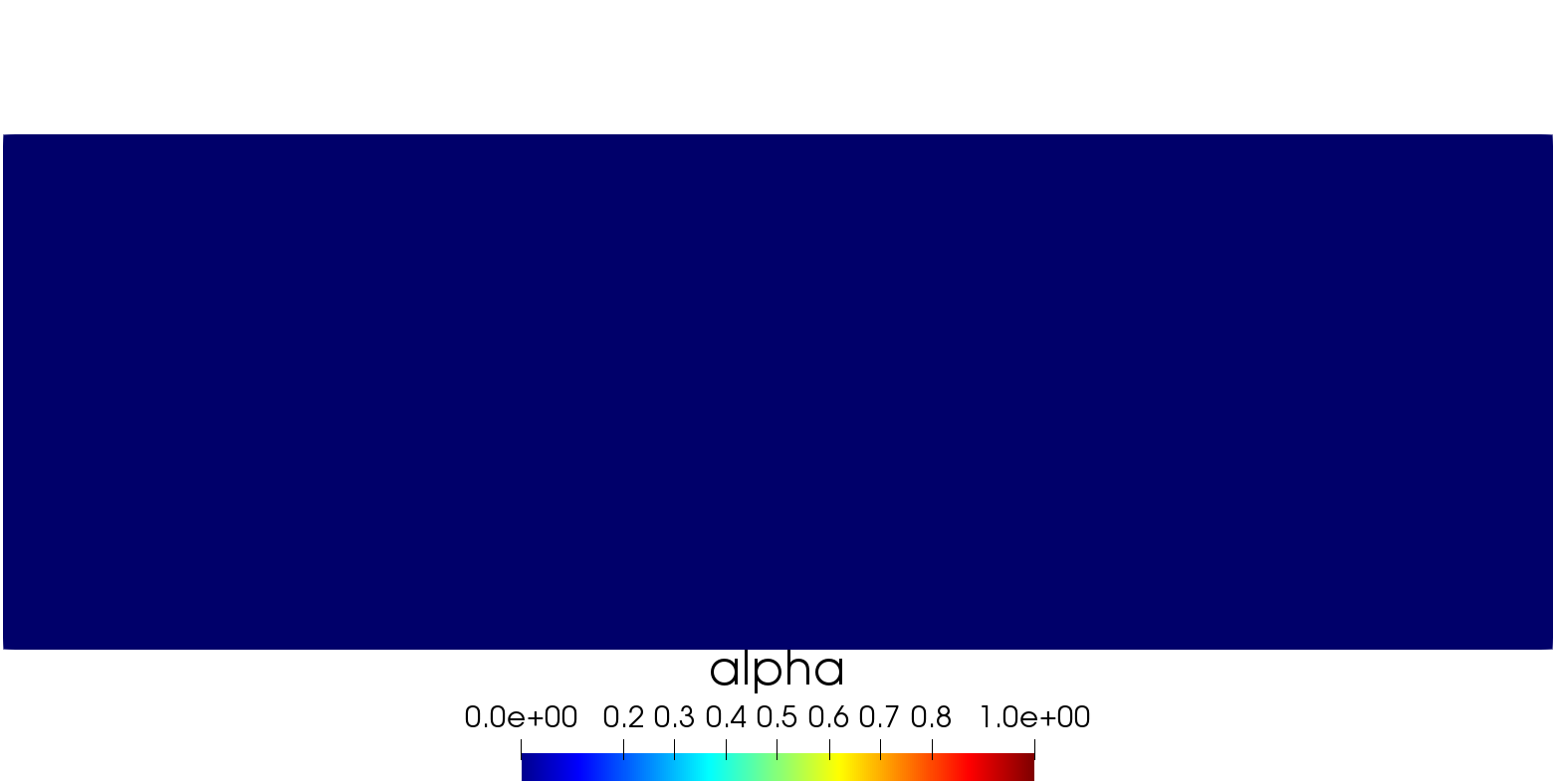}
		\caption{$t_i=0$}
	\end{subfigure}
	\begin{subfigure}[b]{0.45\textwidth}
		\centering
		\includegraphics[width=\textwidth]{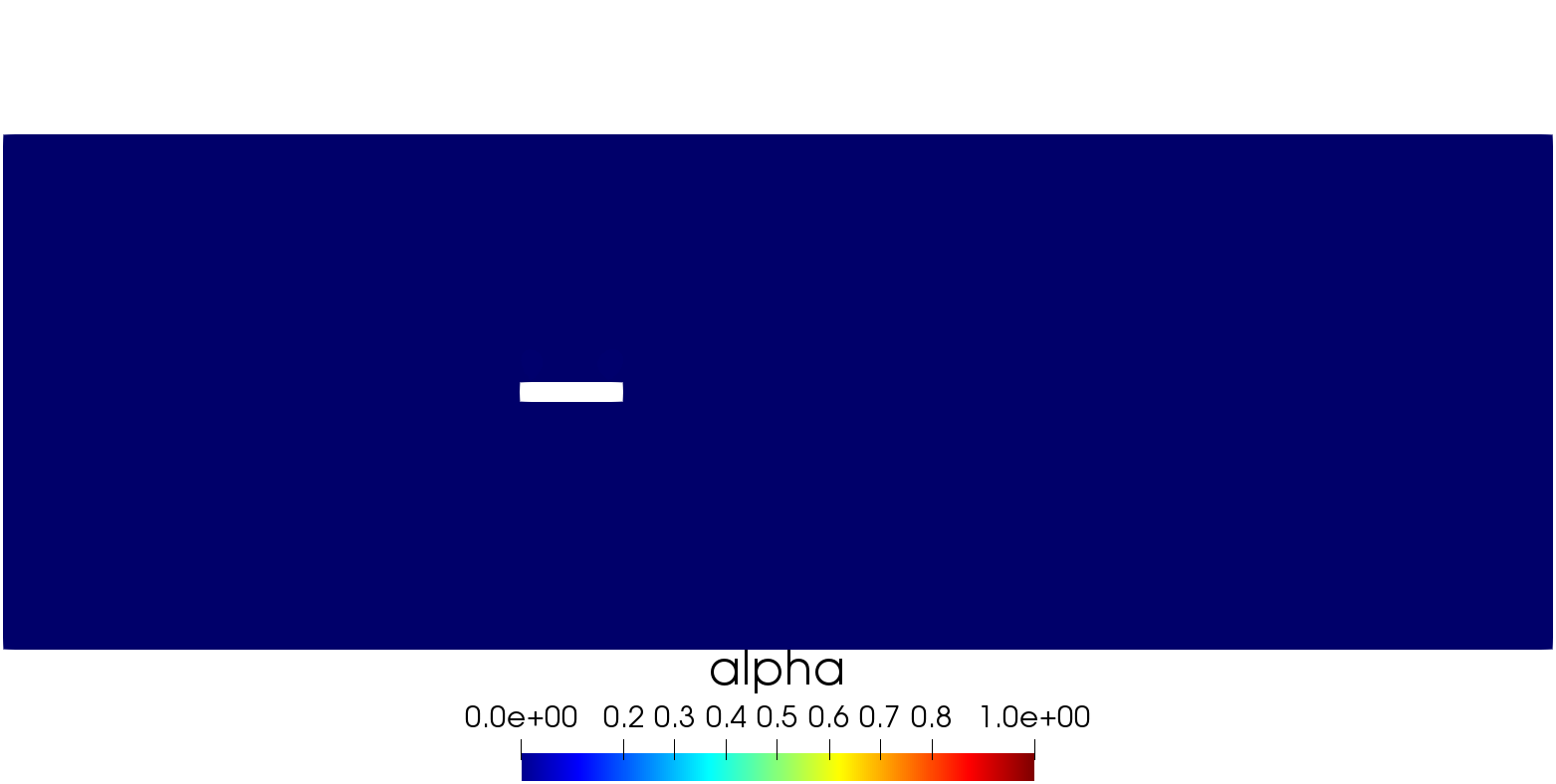}
		\caption{$t_i=5$}
	\end{subfigure}
	\begin{subfigure}[b]{0.45\textwidth}
		\centering
		\includegraphics[width=\textwidth]{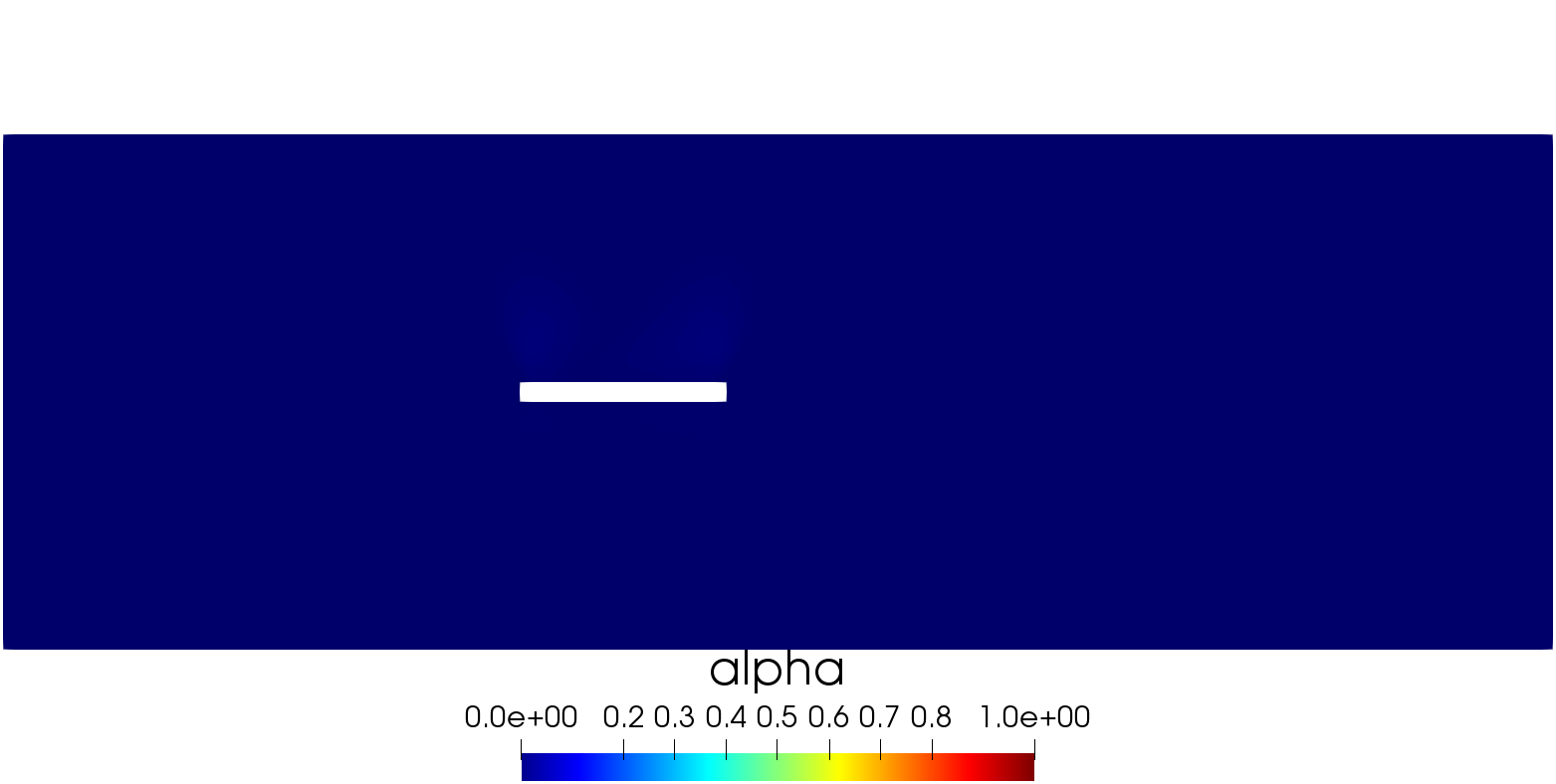}
		\caption{$t_i=10$}
	\end{subfigure}	
	\begin{subfigure}[b]{0.45\textwidth}
		\centering
		\includegraphics[width=\textwidth]{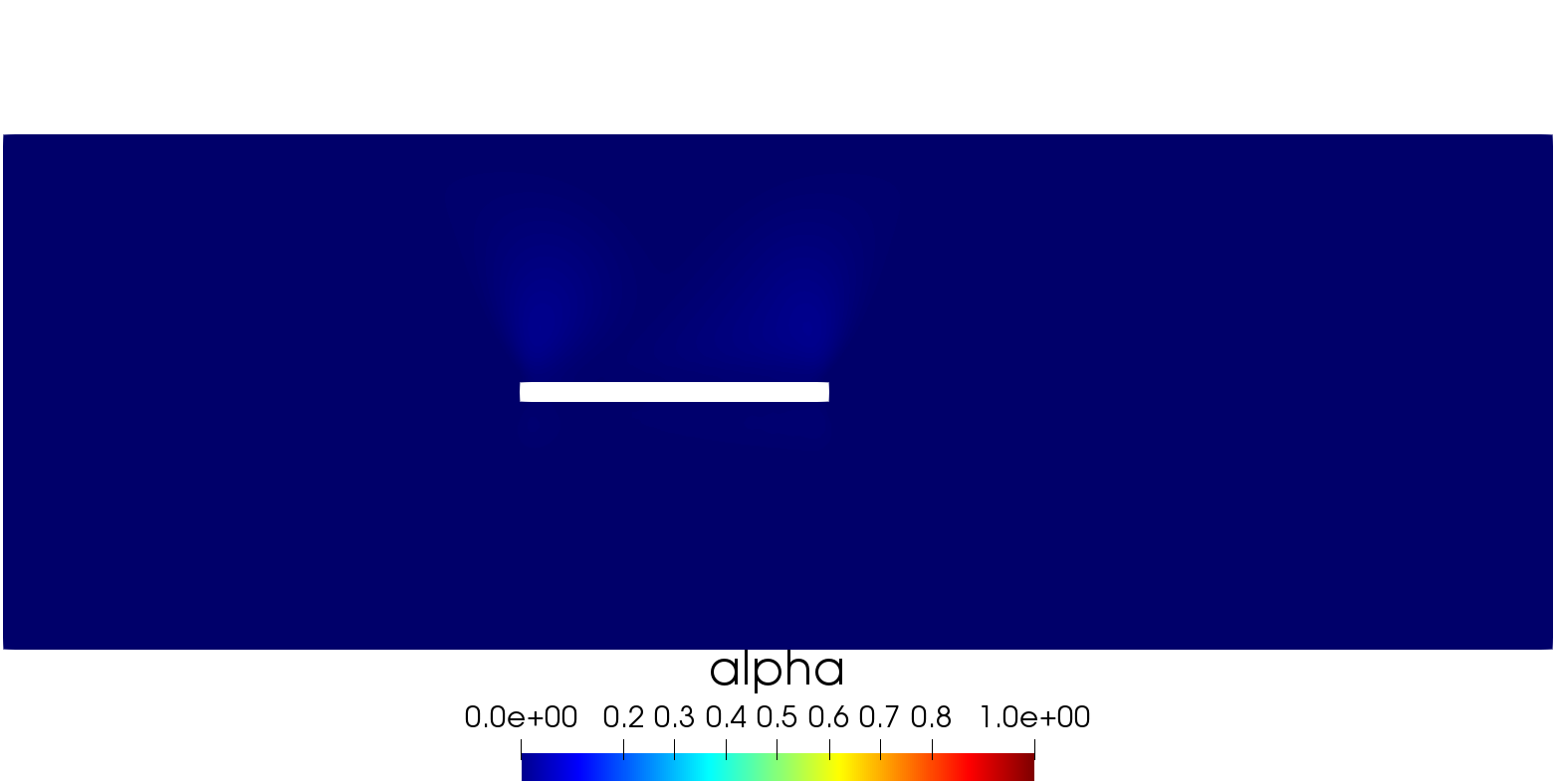}
		\caption{$t_i=15$}
	\end{subfigure}
	\caption{Damage field distribution in the rock mass for $w_1= 10^5\left[\frac{N}{m^3}\right]$}\label{fig5}
\end{figure}

In this test case, we consider our model for the damage problem presented in Section \ref{SCDM}. In Figure \ref{fig5}, the evolution of damage is displayed as the cavity advances in time, for the choice $w_1=10^5\left[\frac{N}{m^3}\right]$. We observe that the level of damage is lower that produced by the gradient damage model, and it is mostly localized around the cavity.

\begin{figure}[h!]
	\centering
	\begin{subfigure}[b]{0.45\textwidth}
		\centering
		\includegraphics[width=\textwidth]{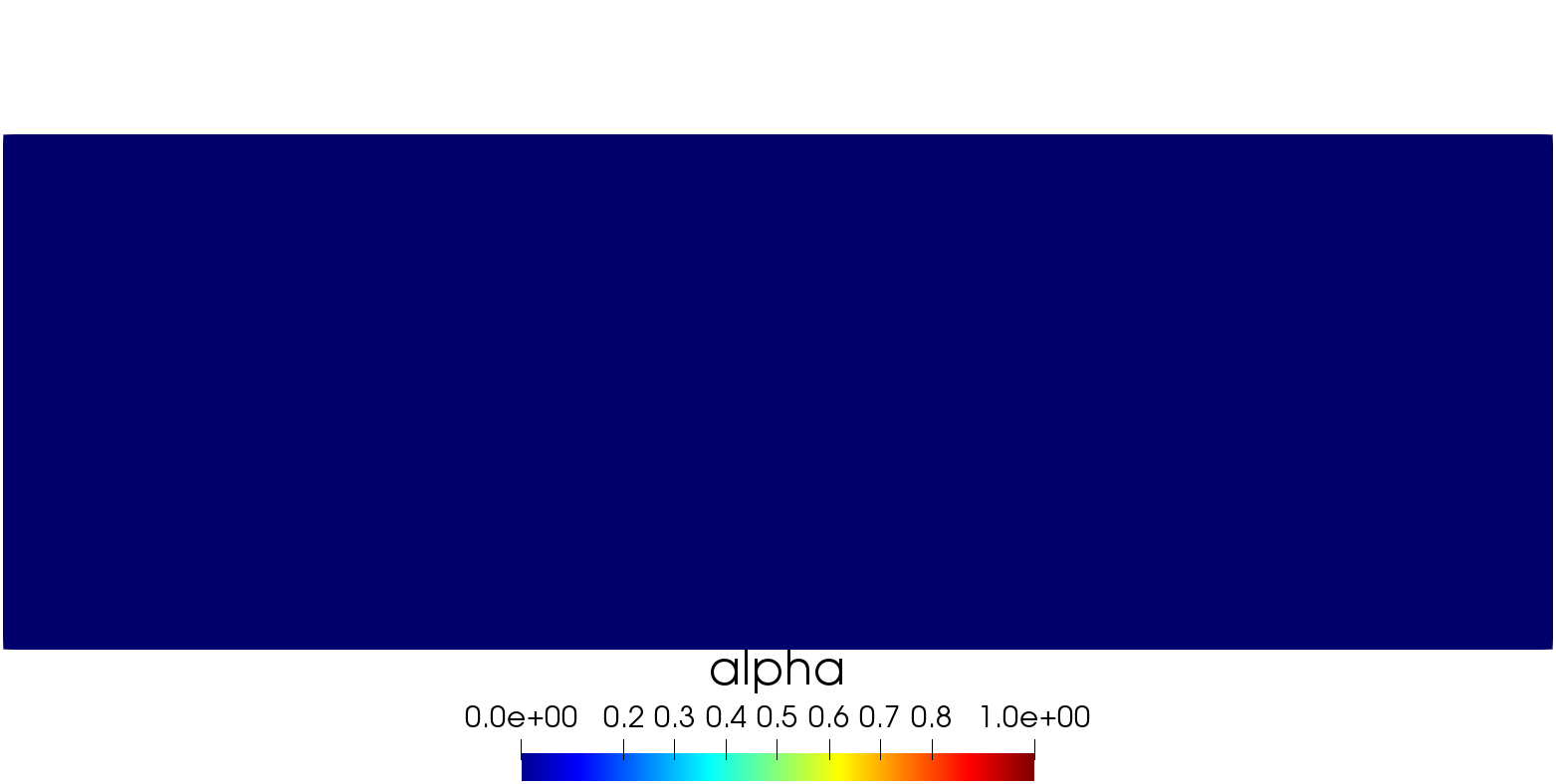}
		\caption{$t_i=0$}
	\end{subfigure}
	\begin{subfigure}[b]{0.45\textwidth}
		\centering
		\includegraphics[width=\textwidth]{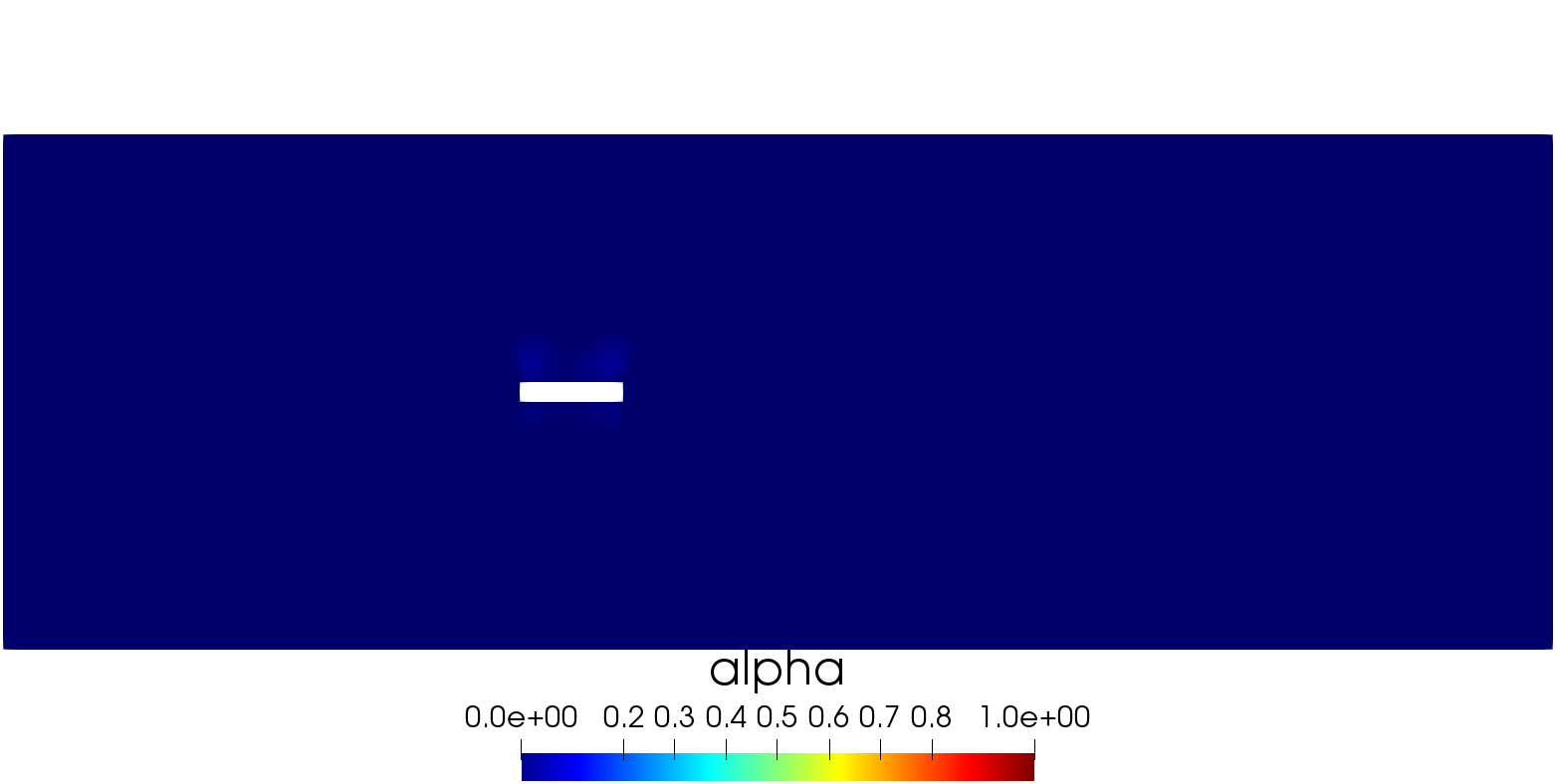}
		\caption{$t_i=5$}
	\end{subfigure}
	\begin{subfigure}[b]{0.45\textwidth}
		\centering
		\includegraphics[width=\textwidth]{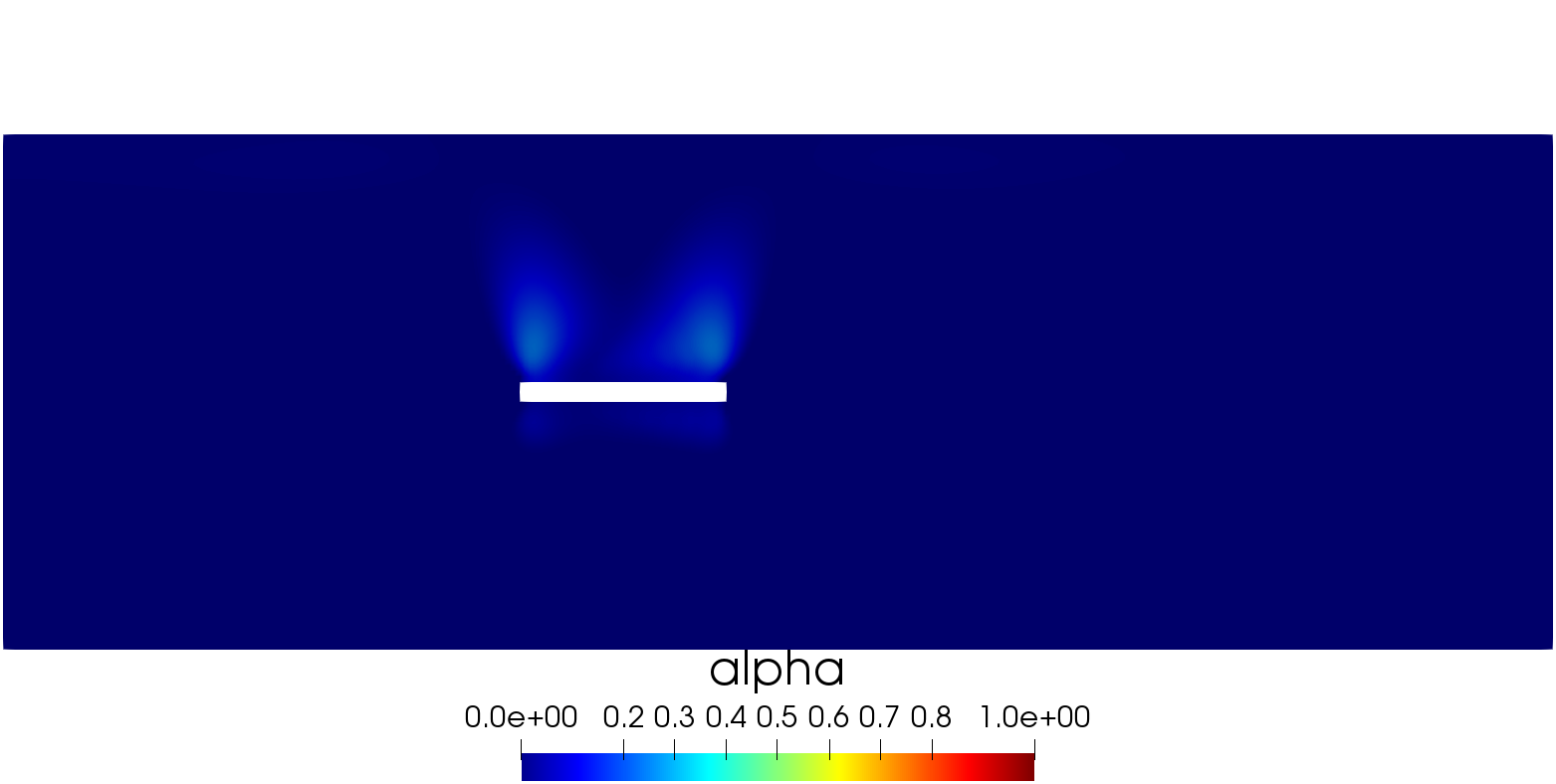}
		\caption{$t_i=10$}
	\end{subfigure}	
	\begin{subfigure}[b]{0.45\textwidth}
		\centering
		\includegraphics[width=\textwidth]{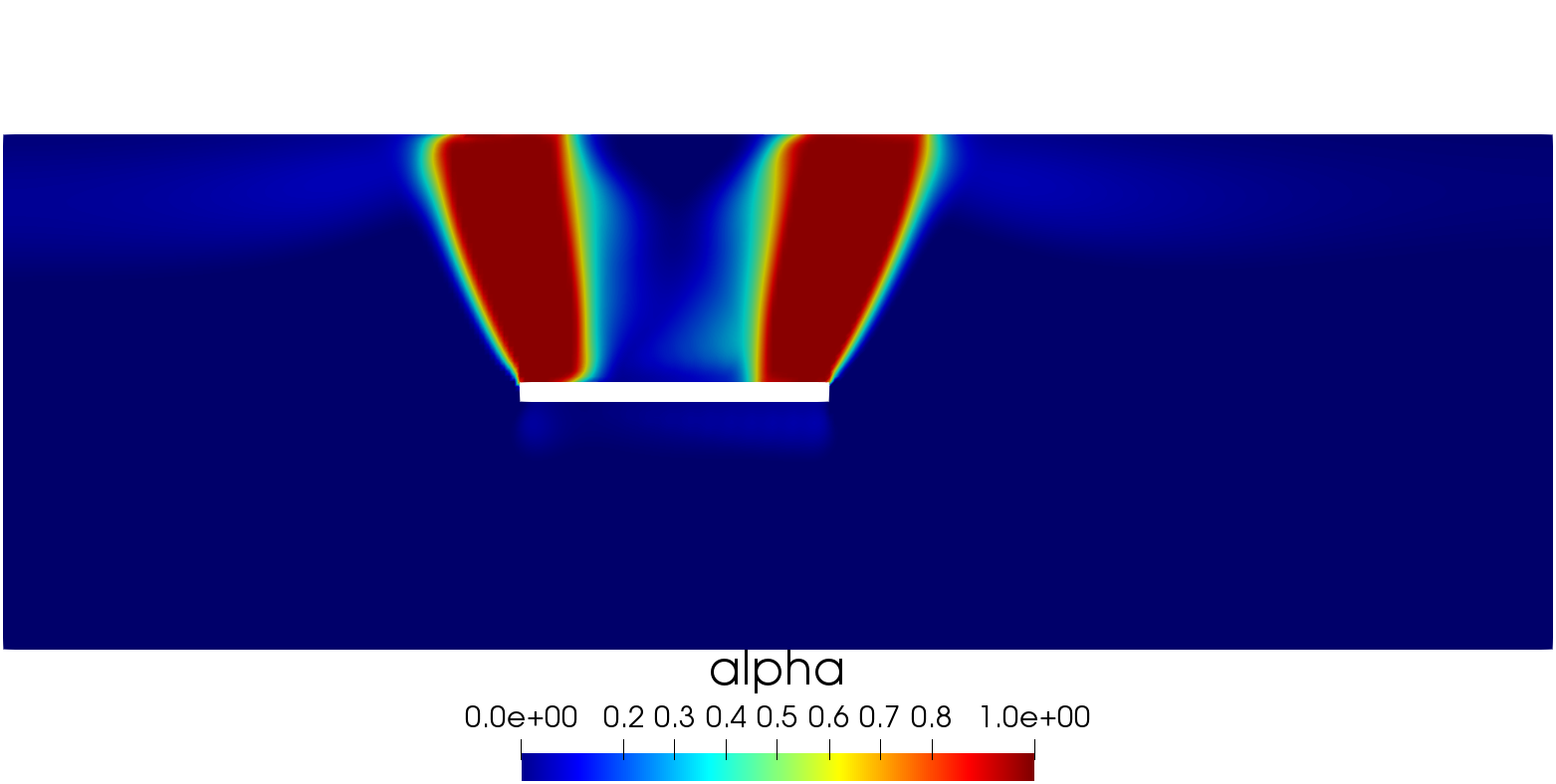}
		\caption{$t_i=15$}
	\end{subfigure}
	\caption{Damage field distribution in the rock mass for $w_1= 10^4\left[\frac{N}{m^3}\right]$.}\label{fig6}
\end{figure}

\begin{figure}[h!]
	\centering
	\begin{subfigure}[b]{0.45\textwidth}
		\centering
		\includegraphics[width=\textwidth]{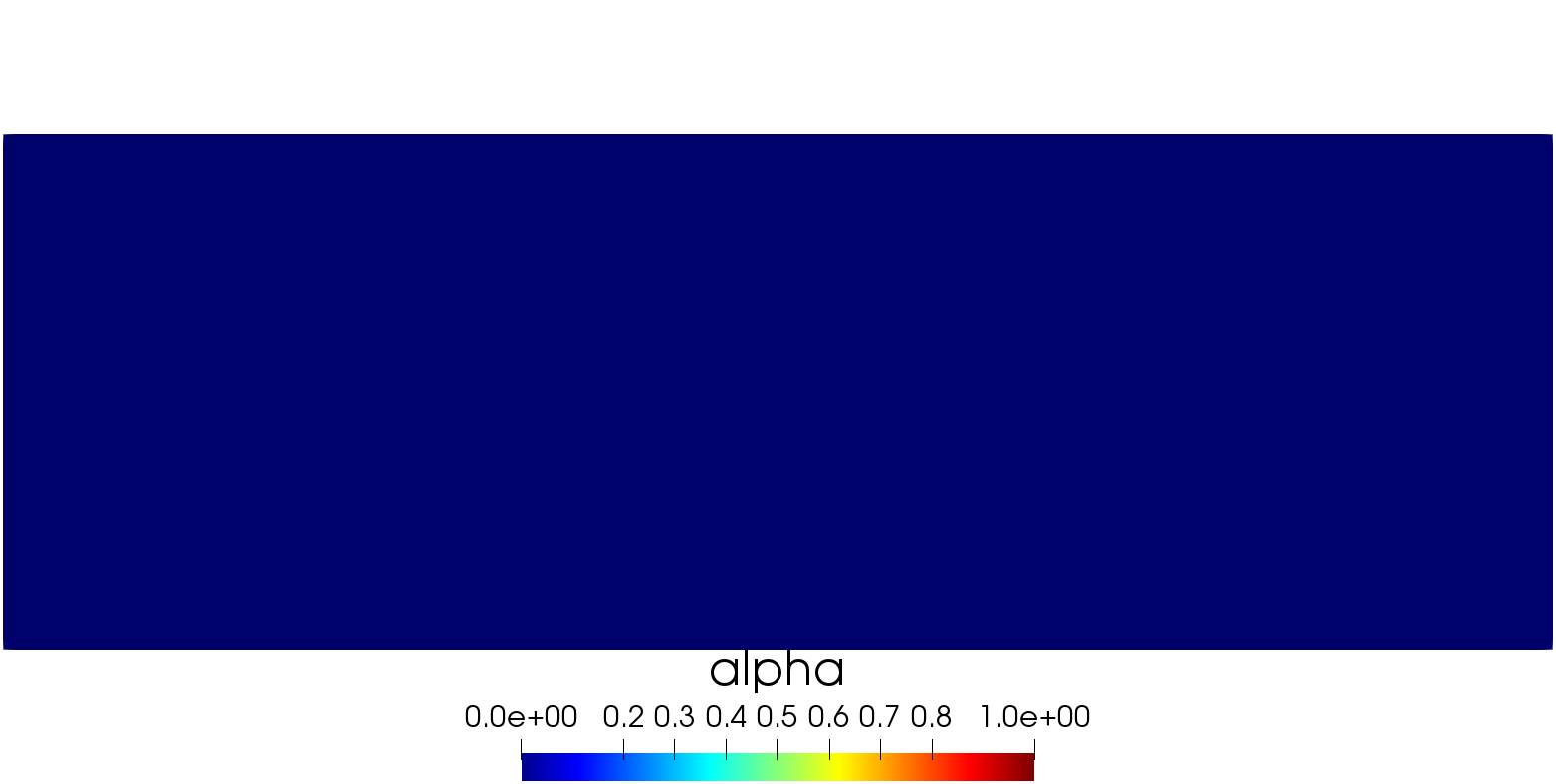}
		\caption{$t_i=0$}
	\end{subfigure}
	\begin{subfigure}[b]{0.45\textwidth}
		\centering
		\includegraphics[width=\textwidth]{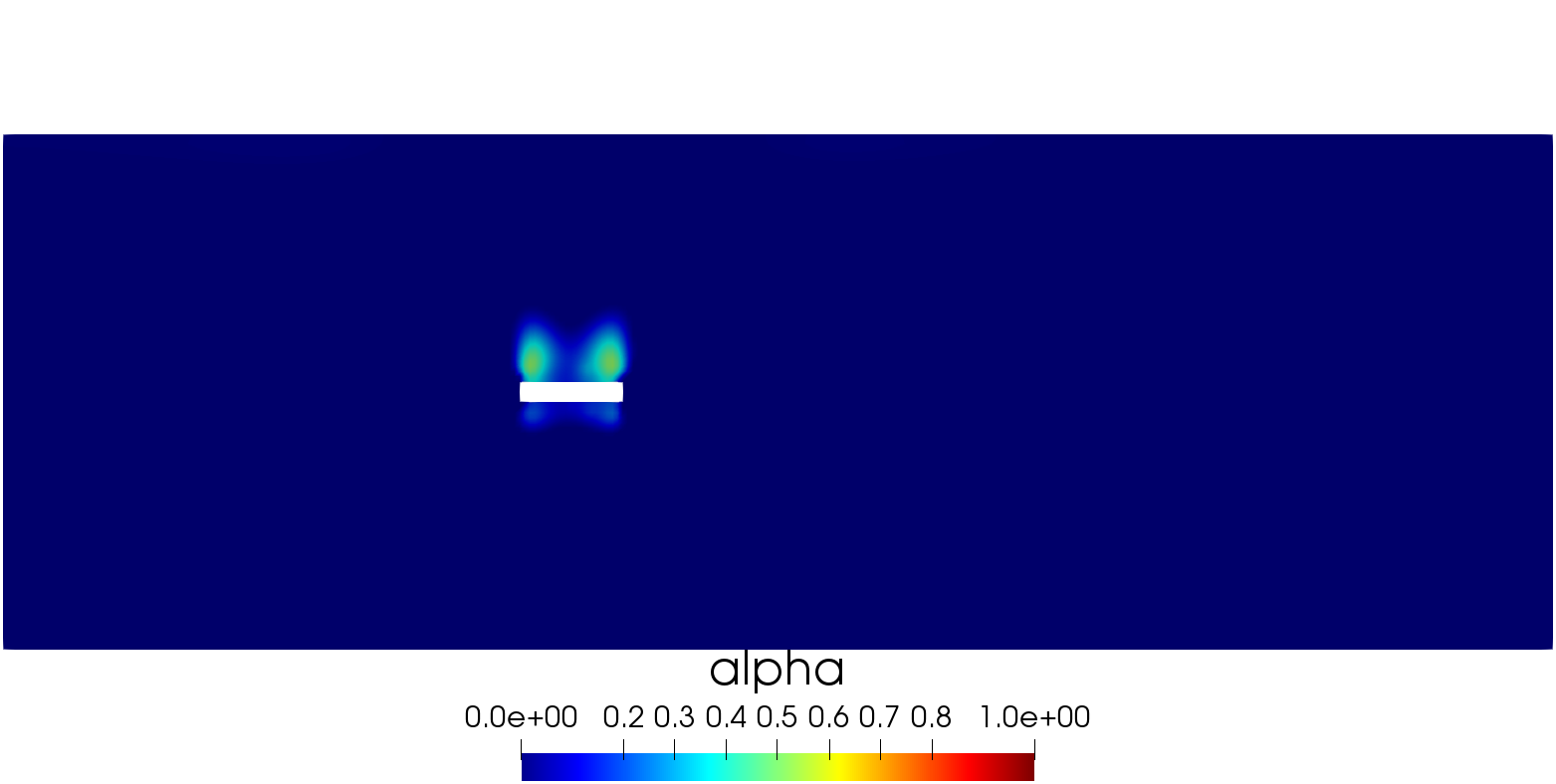}
		\caption{$t_i=5$}
	\end{subfigure}
	\begin{subfigure}[b]{0.45\textwidth}
		\centering
		\includegraphics[width=\textwidth]{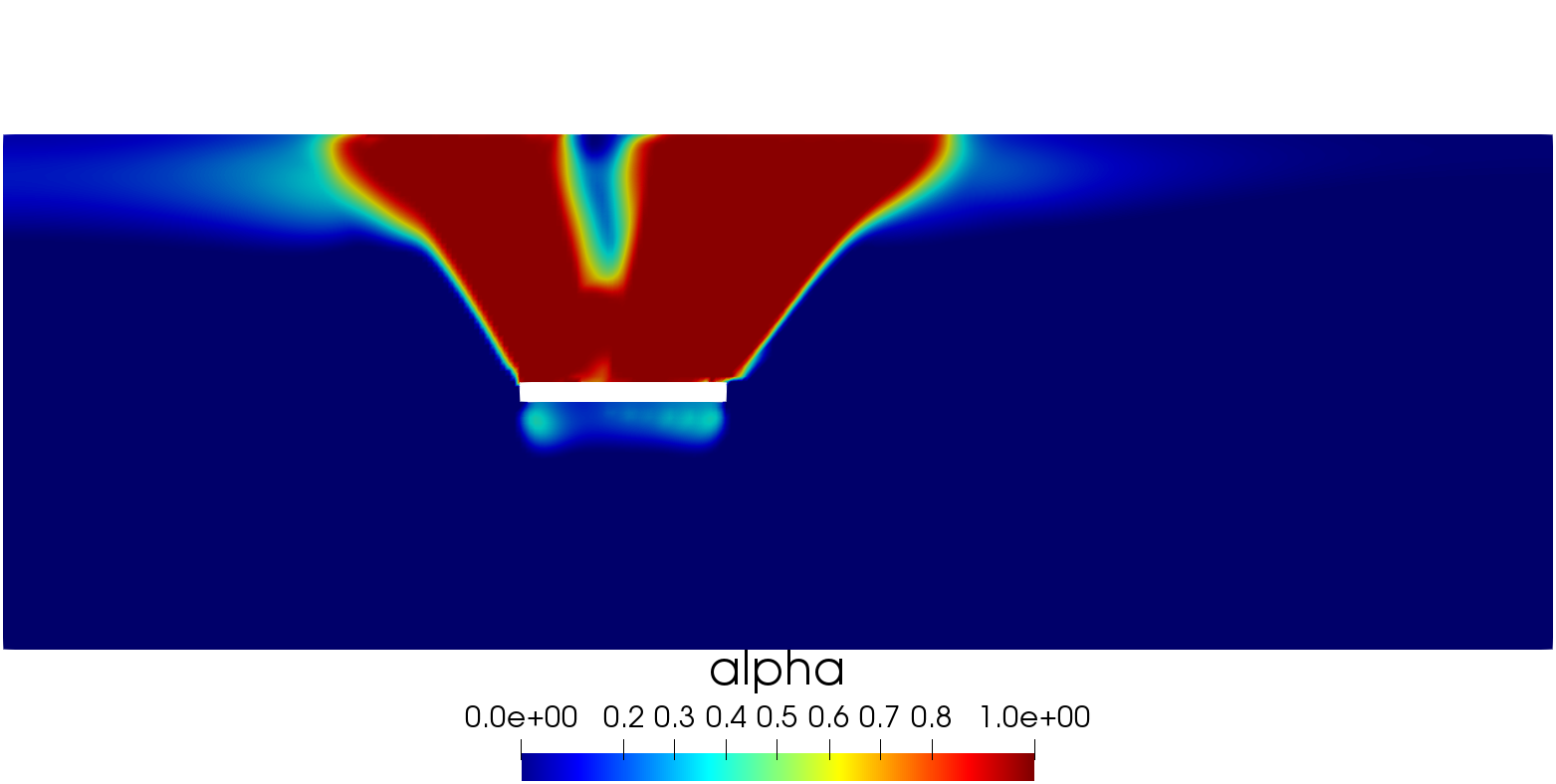}
		\caption{$t_i=10$}
	\end{subfigure}	
	\begin{subfigure}[b]{0.45\textwidth}
		\centering
		\includegraphics[width=\textwidth]{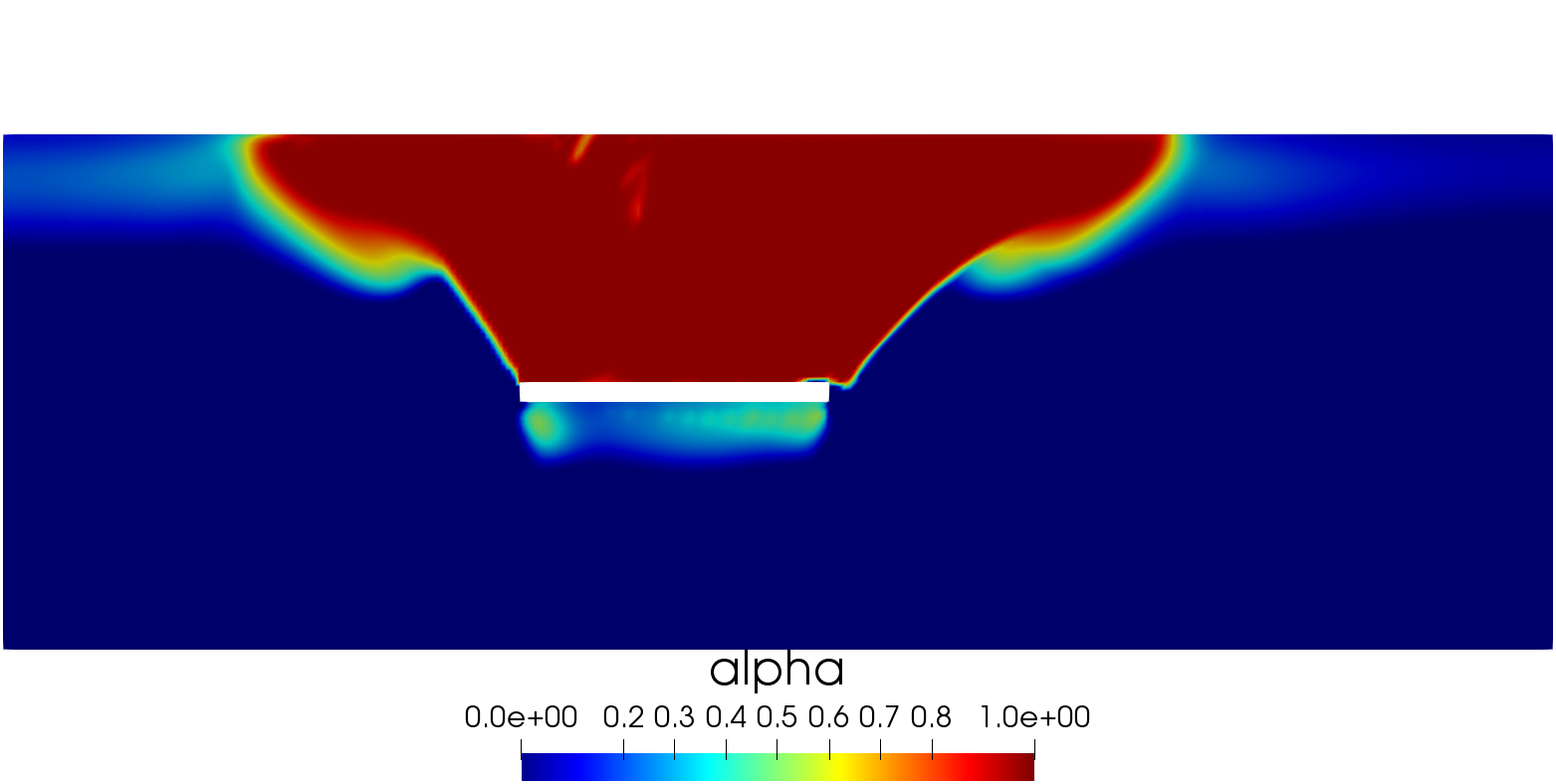}
		\caption{$t_i=15$}
	\end{subfigure}
	\caption{Damage field distribution in the rock mass for $w_1=10^3\left[\frac{N}{m^3}\right]$.}\label{fig7}
\end{figure}

Figures \ref{fig6} and \ref{fig7} show the evolution of damage for $w_1=10^4\left[\frac{N}{m^3}\right]$ and $w_1=10^3\left[\frac{N}{m^3}\right]$ respectively. Damage remains localized around the cavity. As the latter advances in time, damage gets distributed above the ceiling of the cavity, consistently with what is expected in block caving.  

\subsubsection{Sensitivity analysis}

As seen in the damage criterion \eqref{SCdamcrit}, the parameter $\kappa$ plays a fundamental role in the propagation of damage in the rock mass, as it controls the contribution of the spherical and deviatoric part of the stress tensor in the damage criterion.

\begin{figure}[h!]
	\centering
	\begin{subfigure}[b]{0.45\textwidth}
		\centering
		\includegraphics[width=\textwidth]{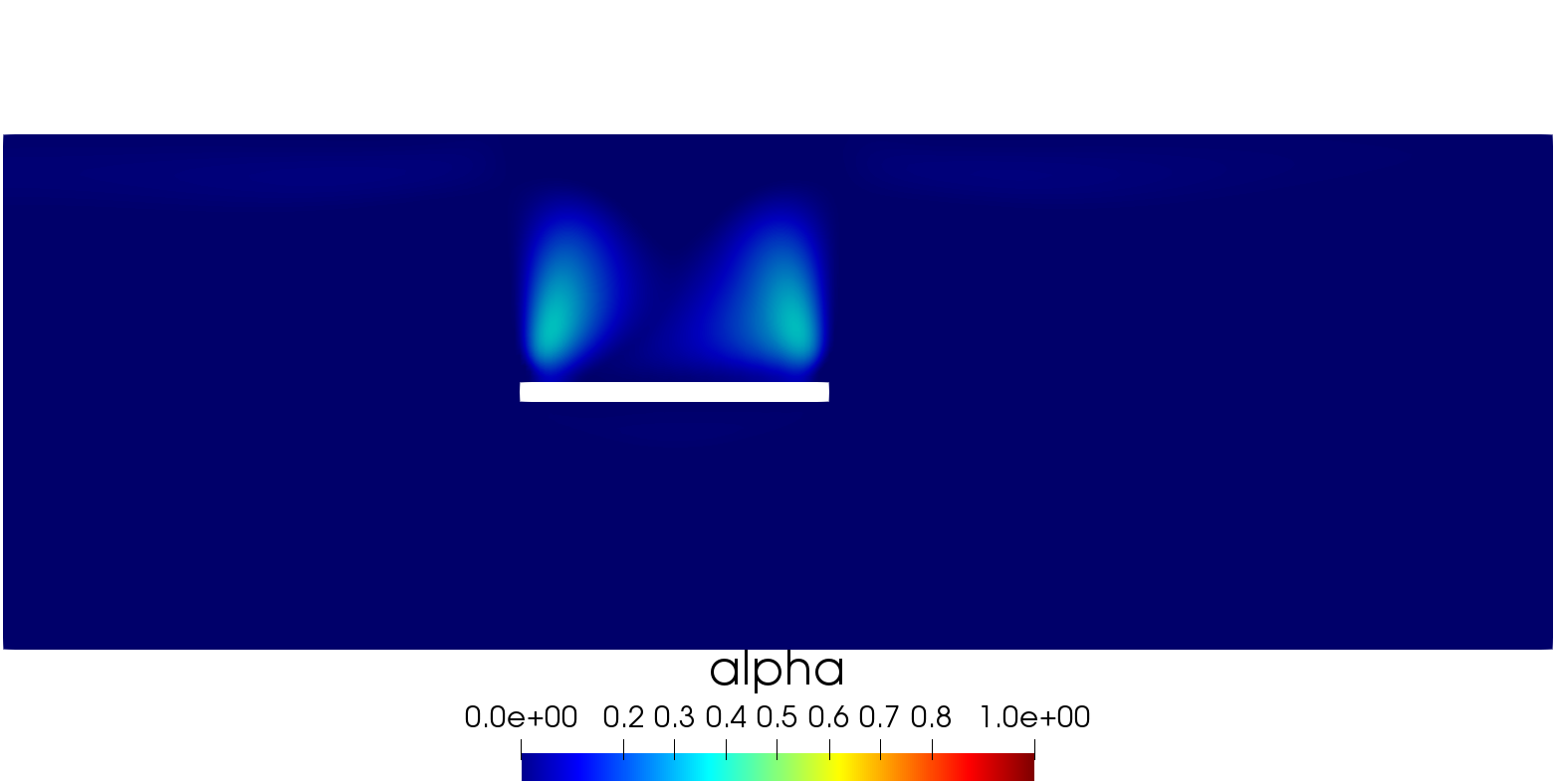}
		\caption{$\kappa=2.0$}
	\end{subfigure}
	\begin{subfigure}[b]{0.45\textwidth}
		\centering
		\includegraphics[width=\textwidth]{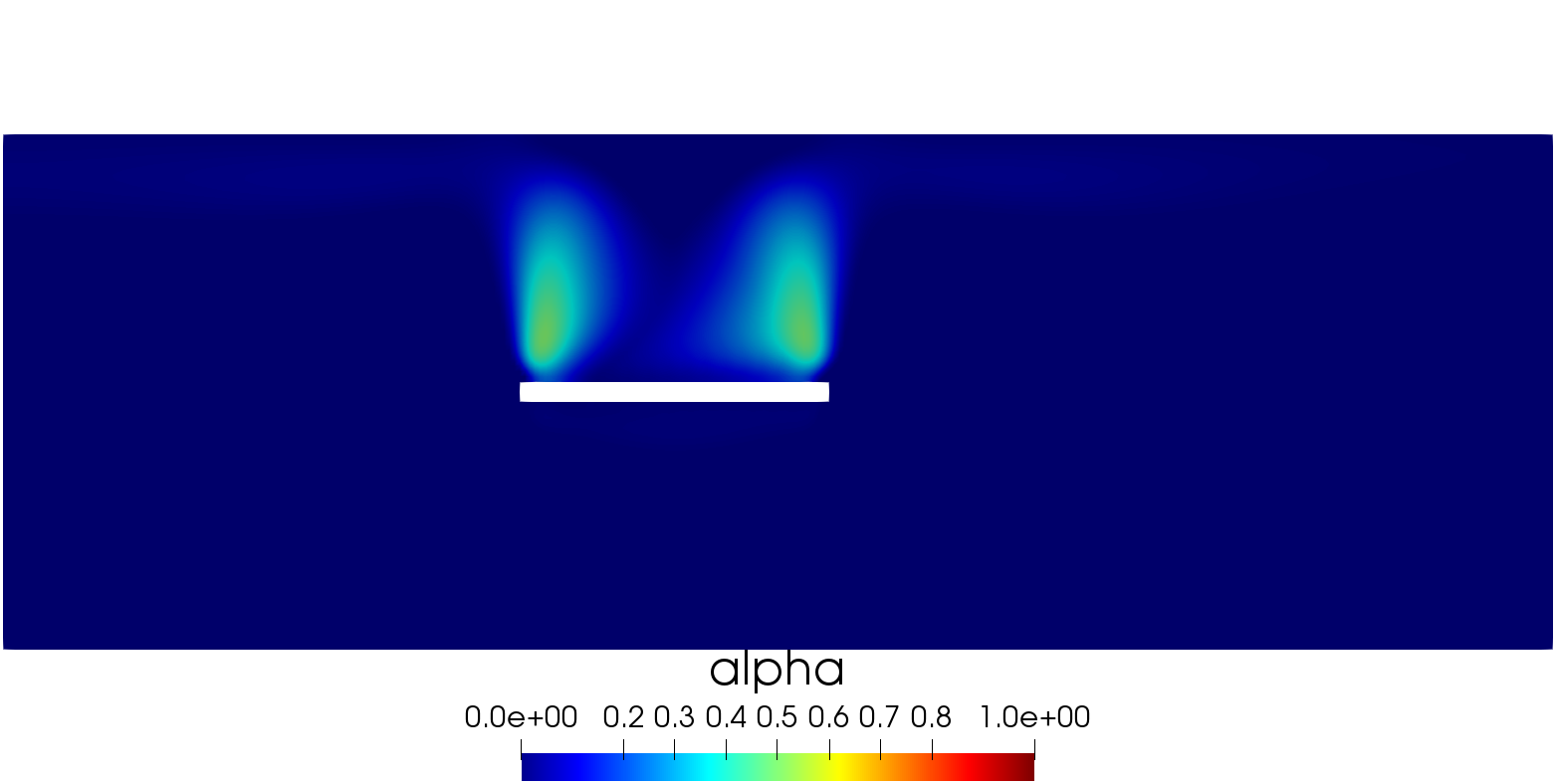}
		\caption{$\kappa=1.5$}
	\end{subfigure}
	\begin{subfigure}[b]{0.45\textwidth}
		\centering
		\includegraphics[width=\textwidth]{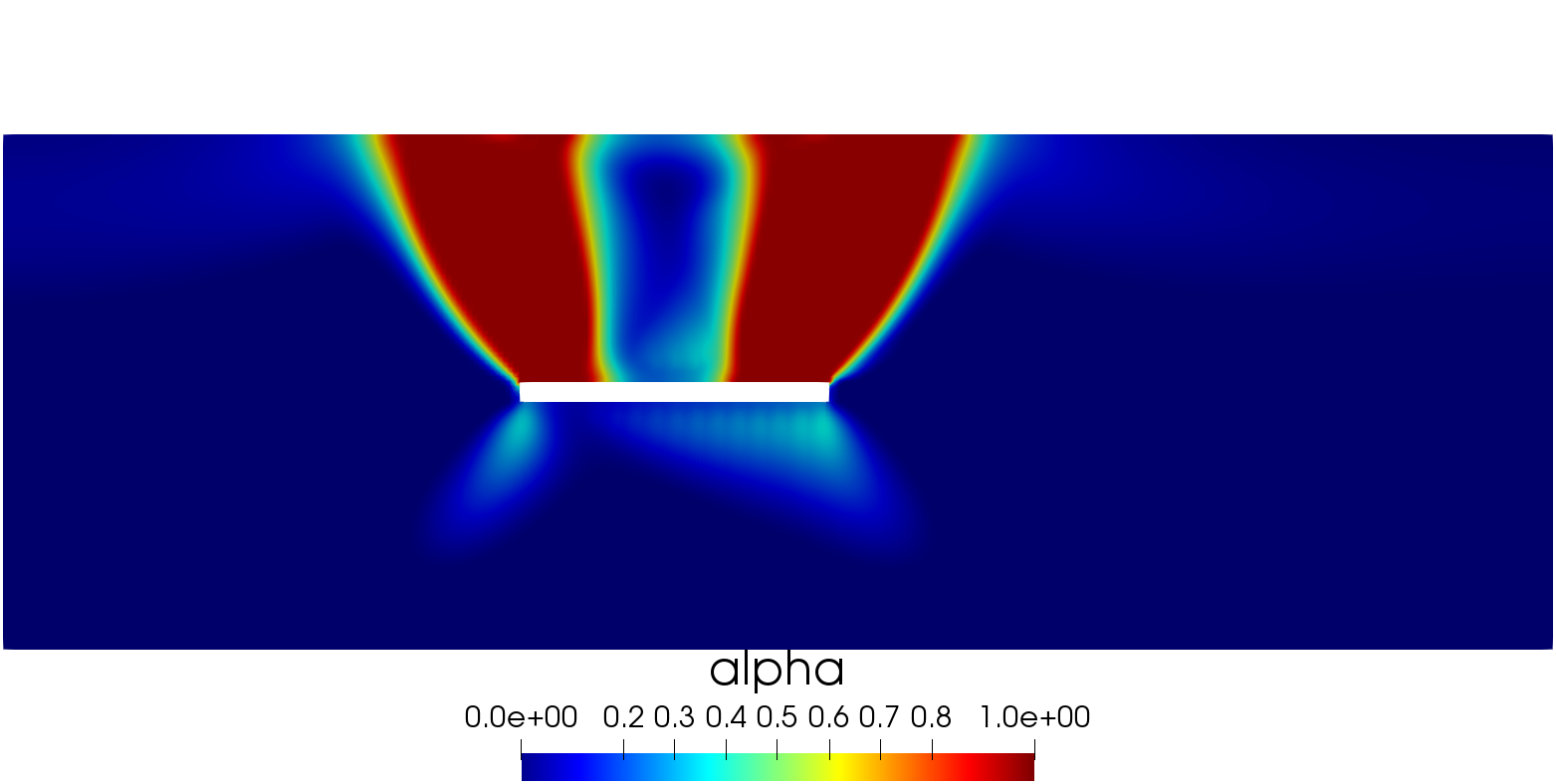}
		\caption{$\kappa=0.5$}
	\end{subfigure}	
	\begin{subfigure}[b]{0.45\textwidth}
		\centering
		\includegraphics[width=\textwidth]{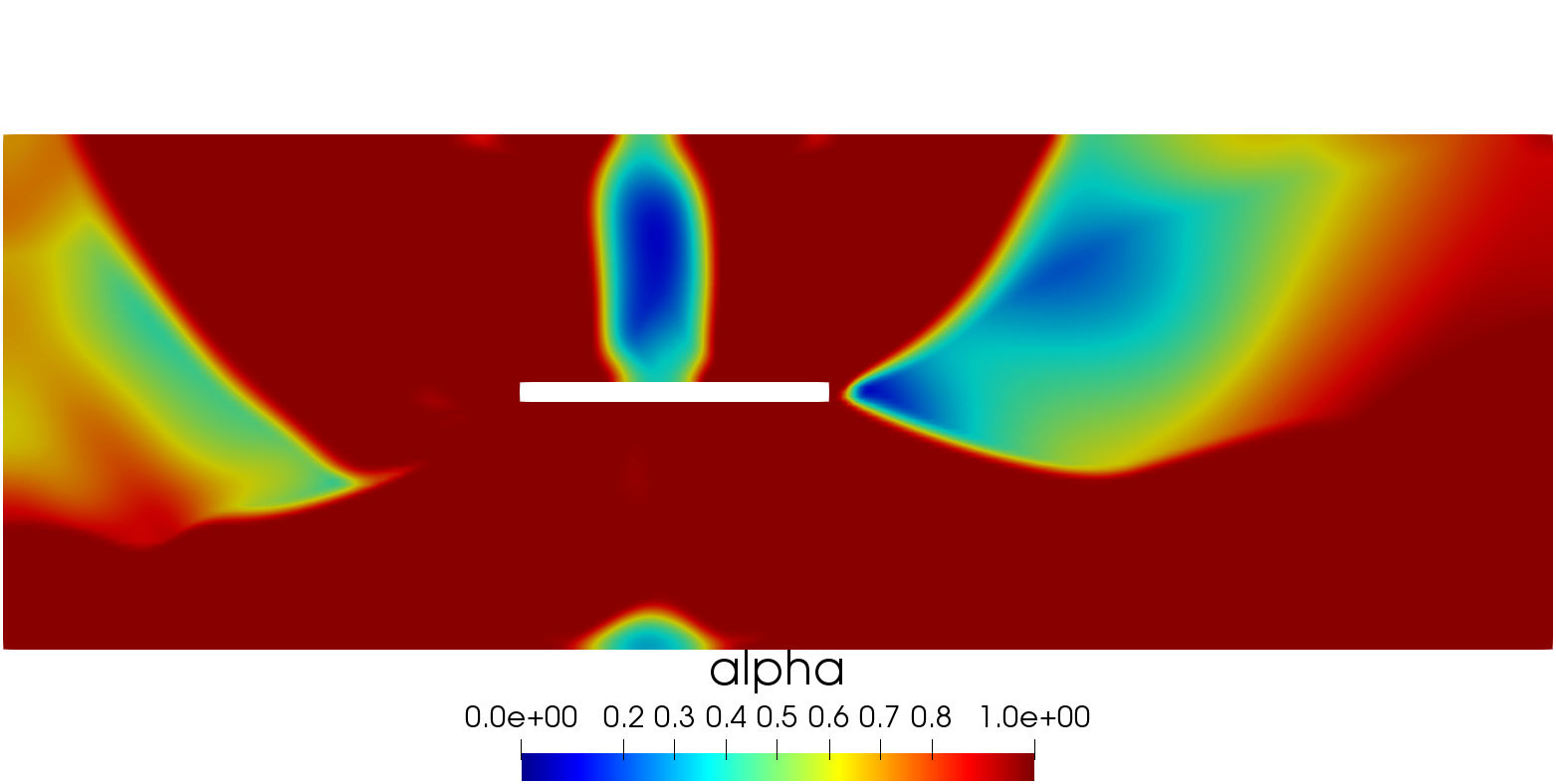}
		\caption{$\kappa=0.2$}
	\end{subfigure}
	\caption{Damage field distribution in the rock mass for $w_1=10^4\left[\frac{N}{m^3}\right]$, with $t_i=15$.}\label{fig8}
\end{figure}

Figure \ref{fig8} displays results for several values of $\kappa$. As this parameter decreases, damage gets distributed throughout the whole domain, starting from the ceiling of the cavity to the upper surface of the rock mass. When $\kappa$ becomes smaller than 1, damage appears in the bottom of the rock mass and invades the whole domain. This is consistent with the fact that for such values of $\kappa$, the damage criterion is less to the deviatoric part of the stresses, whereas the main forcing term is compressive. On the other hand, larger values of $\kappa$ emphasize the deviatoric component of the stress in the damage criterion and privilege shear forces. 

\section{Conclusions}\label{Con}

We studied damage models for representing a process of block caving, based on the gradient damage model of Marigo and Pham \cite{marigo2016overview}. We propose a modified model, that separates the spherical part and the deviatoric part of the stress tensor $\sigma$ in the damage criterion. Our numerical examples show that this model better reproduces the qualitative features expected in block caving: The fracking of the vault of the extraction cavity, as the cavity gets formed. The isotropic gradient damage model in \cite{marigo2016overview} and the gradient damage model for shear fracture in \cite{lancioni2009variational}, do not produce the same qualitative behavior, in the configuration chosen for our numerical examples: The resulting damage mostly occurs below the extraction cavity.

We plan to conduct further investigations so as to build a more realistic model for block caving. These include performing numerical experiment in 3D to account for the effect of edge singularities. We also plan to model the transition for the solid rock to granular phase as the rock mass gets fully damaged, in the spirit of \cite{marigo2019micromechanical}.

\bibliography{mybibfile}
\bibliographystyle{amsplain}

\end{document}